\documentclass{amsart}

\usepackage{amssymb}

\usepackage[all]{xy}

\newtheorem{thm}{Theorem}

\newtheorem{lem}[thm]{Lemma}
\newtheorem{cor}[thm]{Corollary}

\newtheorem{prop}[thm]{Proposition}

\newtheorem{conj}[thm]{Conjecture}
   
\theoremstyle{definition}
\newtheorem{defn}[thm]{Definition}

\newtheorem{say}[thm]{}

\newtheorem{prob}[thm]{Problem}
\newtheorem{ques}[thm]{Question}    

\newtheorem{rem}[thm]{Remark}          

\newtheorem*{ack}{Acknowledgments}      

\newtheorem{defn-thm}[thm]{Definition--Theorem}  
\newtheorem{defn-lem}[thm]{Definition--Lemma}  

\theoremstyle{remark}

\renewcommand{\c}[0]{{\mathbb C}}  

\renewcommand{\o}[0]{{\mathcal O}} 
\newcommand{\z}[0]{{\mathbb Z}}
\newcommand{\n}[0]{{\mathbb N}}
\renewcommand{\r}[0]{{\mathbb R}}

\newcommand{\s}[0]{{\mathbb S}}

\newcommand{\p}[0]{{\mathbb P}}

\newcommand{\q}[0]{{\mathbb Q}}

\newcommand{\qtq}[1]{\quad\mbox{#1}\quad}

\newcommand{\pic}[0]{\operatorname{Pic}}
\newcommand{\gal}[0]{\operatorname{Gal}}

\newcommand{\rank}[0]{\operatorname{rank}}
\newcommand{\mult}[0]{\operatorname{mult}}

\newcommand{\supp}[0]{\operatorname{Supp}}    
\newcommand{\red}[0]{\operatorname{red}}    
\newcommand{\codim}[0]{\operatorname{codim}}    
\newcommand{\im}[0]{\operatorname{im}}    
\newcommand{\proj}[0]{\operatorname{Proj}}

\newcommand{\aut}[0]{\operatorname{Aut}}

\newcommand{\sing}[0]{\operatorname{Sing}}

\newcommand{\res}[0]{\operatorname{\mathcal R}}

\newcommand{\onto}[0]{\twoheadrightarrow}

\newcommand{\tsum}[0]{\textstyle{\sum}}

\def\into{\DOTSB\lhook\joinrel\to}

\def\loccoh#1.#2.#3.#4.{H^{#1}_{#2}(#3,#4)}

\DeclareMathAlphabet{\mathchanc}{OT1}{pzc}%
                                {m}{it}

\newcommand{\SL}{\mathrm{SL}}

\usepackage[all]{xy}\xyoption{dvips}

\def\CC{\mathcal C}

\newcommand{\mdeg}[0]{\operatorname{Deg}}
\newcommand{\dres}[0]{\operatorname{\mathcal D\mathcal R}} 
\newcommand{\tprod}[0]{\textstyle{\prod}} 

\begin{document}
\bibliographystyle{amsalpha}

\title{Links of complex analytic singularities}
\author{J\'anos Koll\'ar}

\maketitle

Let $X$ be a complex algebraic or analytic variety.
Its local topology near  a point $x\in X$
is completely described by  its {\it link}  $L(x\in X)$,
which is obtained as  the intersection of
$X$ with a sphere of radius $0<\epsilon\ll 1$ centered at $x$.
The intersection of $X$ with the closed ball  of radius $\epsilon$ 
centered at $x$ is 
homeomorphic to the cone over $L(x\in X)$; cf.\ \cite[p.41]{gm-book}.

If $x\in X$ is a smooth point then its link is a
sphere of dimension $2\dim_{\c} X-1$. Conversely, if $X$ is a normal surface and
 $L(x\in X)$ is a sphere then $x$ is a smooth point \cite{mumf-top},
but this fails in higher dimensions \cite{briesk-exmps}.

The aim of this survey  is to study in some sense
the opposite question: we are interested in the ``most complicated'' links. 
 In its  general form, the question is the following.

\begin{prob}\label{general.prob}
 Which topological spaces can be links of 
 complex algebraic or analytic singularities?
\end{prob}

If $\dim X=1$, then the possible links are  disjoint unions
of circles.  The answer is much more complicated in higher dimensions
and we focus on isolated singularities from now on, though many results hold
for non-isolated  singularities as well. Thus the link
 $L(x\in X)$ is a (differentiable) manifold of (real)
 dimension $2\dim_{\c} X-1$.

Among the simplest singularities are the cones over smooth
projective varieties. Let $Z\subset \p^N$ be a smooth
projective variety and $X:=\operatorname{Cone}(Z)\subset \c^{N+1}$ 
the cone over $Z$
with vertex at the origin. Then
$L(0\in X)$ is a circle bundle over $Z$ whose first Chern class
is the hyperplane class. Thus the link of the vertex of  
$\operatorname{Cone}(Z) $
is completely described by the base $Z$ and by the 
hyperplane class $[H]\in H^2(Z, \z)$. 

Note that a singularity  $0\in X\subset \c^N$ is a cone
iff it can be defined by homogeneous equations.
One gets a much larger class of singularities if we consider
homogeneous equations where different variables  have different 
degree (or weight). 

For a long time it was believed that links of isolated singularities 
are ``very similar'' to links of cones and weighted cones.
 The best illustration of this
is given by the complete description of  links of surface
singularities given in \cite{neumann-pl}. 
Cones give circle bundles over Riemann surfaces
and weighted cones give Seifert bundles
over  Riemann surfaces.
General links are more complicated but they are all obtained
by gluing  Seifert bundles over Riemann surfaces with boundary.
 These are definitely more  complicated
than  Seifert bundles, but much simpler than
general 3--manifolds. In particular, hyperbolic 3--manifolds
-- which comprise the largest and most complicated class --
do not occur as links.

Important examples of the similarity of general links to smooth
projective varieties are given by the local Lefschetz theorems,
initiated by Grothendieck \cite{sga2} and developed much further
subsequently; see \cite{gm-book} for a detailed treatment.

As another illustration,
the weights of the mixed Hodge structure on the cohomology groups of links
also follow the same pattern for general links as for links of cones, see
\cite{MR928298} or \cite[Sec.6.3]{PetersSteenbrinkBook}.

These and many other examples led to a viewpoint that was
best summarized in
 \cite[p.26]{gm-book}:
{\it ``Philosophically, any statement about the projective variety 
or its embedding
 really comes from a statement about the singularity at the point of the cone. 
Theorems about projective varieties should be consequences of more general 
theorems about singularities which are no longer required to be conical.''}

Recently this belief was called into question
by \cite{k-fg} which proved that fundamental groups of general links
are very different from fundamental groups of links of cones.
The aim of this paper
is to summarize the  results, present several new theorems
 and review the problems that arise.

Philosophically, the main long term question is to understand the limits of the
above principle. We know that it fails for the fundamental group
but it seems to apply to cohomology groups. It is  unclear if
it applies to simply connected links or not.

The new results rely on a 
 method, considered in \cite{k-t-lc-exmp},
 to construct singularities using their resolution.
 By Hironaka's resolution theorem,
for every isolated singularity $(x\in X)$ there is a
proper, birational morphism
$f:Y\to X$ such that $E:=f^{-1}(x)$ is a simple normal crossing divisor
and $Y\setminus E\to X\setminus \{x\}$ is an isomorphism.
The method essentially reverses the resolution process.
That is, we start with a (usually reducible) 
simple normal crossing variety $E$,
embed $E$ into a smooth variety $Y$ and then contract
$E\subset Y$ to a point to obtain  $(x\in X)$.
If $E$ is  smooth, this is essentially the cone construction.

This approach has been one of the standard ways to construct
surface singularities but it has not been investigated in higher dimensions
until recently. There were probably two reason for this.
First, if $\dim X\geq 3$ then there is no ``optimal'' choice for
the resolution $f:Y\to X$. Thus the exceptional set $E=f^{-1}(x)$
depends on many arbitrary choices and  it is not easy to extract any
invariant of the singularity from $E$;
see, however,  Definition \ref{dres.defn}. 
Thus any construction starting with $E$
seemed rather arbitrary.
Second, the above philosophy suggested that one should not get
anything substantially new this way.

The first indication that this method is worth exploring
was given in \cite{k-t-lc-exmp} where it was used to construct
new examples of terminal and log canonical singularities that
contradicted earlier expectations.

A much more significant application was given in \cite{k-fg}.
Since in higher dimensions a full answer to Problem 
\ref{general.prob} may well be impossible to give, 
it is sensible to focus  on some special aspects.
A very interesting question turned out to be the following.

\begin{prob}\label{general.fg.prob}
 Which groups occur as fundamental groups of links of 
 complex algebraic or analytic singularities?
\end{prob}

Note that the fundamental groups of smooth projective varieties
are rather special; see \cite{abckt} for a  survey. 
 Even the fundamental groups of smooth  quasi projective
varieties are  quite restricted \cite{Morgan, kap-mil, Corlette-Simpson, DPS}.
By contrast fundamental groups of links are arbitrary.

\begin{thm}\cite{k-fg} \label{link.thm}
For every finitely presented group $G$ there is an
isolated,  complex singularity $\bigl(0\in X_{G}\bigr)$
with link $L_{G}$
such that  $\pi_1\bigl(L_{G}\bigr)\cong G$.
\end{thm}

Note that once such a singularity exists, a local Lefschetz--type
theorem (cf.\  \cite[Sec.II.1.2]{gm-book})
 implies that the link of a general 3-dimensi\-onal
hyperplane section has the same fundamental group. 
\medskip

There are two natural directions to further develop this result:
one can connect properties of the fundamental group of a link
to algebraic or analytic properties of a singularity
and one can investigate further the topology of the links
or of the resolutions.

In the first direction, the following result
 answers a question of Wahl.

\begin{thm} \label{cm.main.thm.intro}
 For a finitely presented group $G$ the following are equivalent.
\begin{enumerate}
\item $G$ is $\q$-perfect, that is, 
its largest abelian quotient is finite.
\item $G$ is the fundamental group of the link of an isolated 
Cohen--Macaulay singularity (\ref{CM.defn}) of dimension $\geq 3$.
\end{enumerate}
\end{thm}

One can study the local topology of
$X$ by  choosing a resolution of
 singularities  $\pi:Y\to X$ such that
$E_x:=\pi^{-1}(x)\subset Y$ is a simple normal crossing divisor
and then relating the topology of $E_x$ to the topology of the
 link  $L(x\in X)$.

The topology of a simple normal crossing divisor
$E$ can in turn be understood in 2 steps.
First, the $E_i$ are smooth projective varieties,
and their topology is much studied.
A second layer of complexity comes from how the
components $E_i$ are glued together. This gluing process can be
naturally encoded by a finite cell  complex  ${\mathcal D}(E)$,
called the 
{\it dual  complex} or  {\it dual graph} of $E$.

\begin{defn}[Dual complex]\label{snc-nc.sing.defn}
Let $E$ be a variety
with irreducible components $\{E_i: i\in I\}$.
We say that $E$ is a {\it simple normal crossing} variety
(abbreviated as {\it snc}) if the $E_i$ are smooth and 
 every point $p\in E$ has an open (Euclidean) neighborhood $p\in U_p\subset E$
and an  embedding $U_p\into \c^{n+1}$ 
such that  the image of $U_p$ is an open subset of
the union of coordinate hyperplanes $ (z_1\cdots z_{n+1}=0)$.
A {\it stratum}
of $E$ is any irreducible component of an intersection
$\cap_{i\in J} E_i$ for some $J\subset I$.

The combinatorics of $E$ is encoded by a cell complex ${\mathcal D}(E)$
 whose vertices are labeled by the
 irreducible components of $E$ and for every stratum
$W\subset \cap_{i\in J} E_i$ we attach a $(|J|-1)$-dimensional cell.
Note that for any $j\in J$ there is a unique  irreducible component
of $\cap_{i\in J\setminus\{j\}} E_i$ that contains $W$;
this specifies the attaching map.
${\mathcal D}(E)$ is called the
{\it dual complex} or  {\it  dual graph} of $E$. 
(Although ${\mathcal D}(E)$ is not a simplicial complex in general, it is an
unordered $\Delta$-complex in the terminology of
\cite[p.534]{hatcher}.)
\end{defn}

\begin{defn}[Dual complexes associated to a singularity]\label{dres.defn}
 Let $X$ be a normal variety and
$x\in X$ a point. Choose a resolution of
 singularities  $\pi:Y\to X$ such that
$E_x:=\pi^{-1}(x)\subset Y$ is a simple normal crossing divisor.
Thus it has a dual  complex  ${\mathcal D}(E_x)$.

The dual graph of a normal surface singularity has a long history. Higher
dimensional versions appear in 
\cite{MR0506296, MR0466149, MR576865, friedman-etal}
but systematic investigations were started only recently; see
\cite{MR2320738, MR2399025, payne09, payne11}.

It is proved in \cite{MR2320738, MR2399025, 2011arXiv1102.4370A}
 that  the homotopy type of 
${\mathcal D}(E_x)$
is  independent of the resolution $Y\to X$. We denote it by
$\dres(x\in X)$. 
\end{defn}

The proof of  Theorem \ref{link.thm} 
gives singularities for which
the fundamental group of the link is isomorphic to
the fundamental group of $\dres(x\in X)$.
In general, it seems  easier to study $\dres(x\in X)$
than the link and the next theorem  shows that
not just the fundamental group but the whole
homotopy type of 
 $\dres(0\in X)$ can be arbitrary. The additional properties
(\ref{main.thm.1}.2--3) follow from the construction as in
\cite{k-t-lc-exmp, k-fg}.

\begin{thm}\label{main.thm.1}
 Let $T$ be a  connected,  finite  cell complex.
Then there is a normal singularity $(0\in X)$ such that
\begin{enumerate}
\item the  complex $\dres(0\in X)$ is homotopy equivalent to $T$,
\item $\pi_1\bigl(L(0\in X)\bigr)\cong \pi_1(T)$ and
\item if $\pi:Y\to X$ is any resolution then
$R^i\pi_*\o_{Y}\cong H^i(T, \c)$ for $i>0$.
\end{enumerate}
\end{thm}

The fundamental groups of the dual complexes of 
rational singularities (\ref{rtl.defn})
were determined in \cite[Thm.42]{k-fg}.
The next result extends this by determining
 the possible homotopy types of 
 $\dres(0\in X) $.

\begin{thm}\label{main.thm.2}
 Let $T$ be a  connected, finite  cell complex.
Then there is a rational singularity $(0\in X)$
whose dual  complex $\dres(0\in X)$ is homotopy equivalent to $T$
iff $T$ is  $\q$-acyclic, that is, $H^i(T,\q)=0$ for $i>0$.
\end{thm}

As noted in \cite{payne11}, the  dual  complex  $\dres(0\in X)$
can be defined even up-to simple-homotopy equivalence \cite{MR0362320}.
The  proofs given in \cite{k-fg} use Theorem \ref{cor.21.thm}, which in turn
relies on some general theorems of 
\cite{Cairns, Hirsch} that do not seem to give simple-homotopy 
equivalence.\footnote{This problem is settled in \cite{k-dual2}.}

\subsection*{Content of the Sections}{\ }

Cones, weighted cones and the topology of the corresponding links
are discussed in Section \ref{sec.1}.

The plan for the construction of singularities from their resolutions
is outlined in Section \ref{sec.2} and
the rest of the paper essentially fleshes out the details.

In Section \ref{sec.vor} we show that every finite cell complex
is homotopy equivalent to a  Voronoi complex.
These Voronoi complexes are then used to construct
simple normal crossing varieties in Section \ref{sec.snc}.

The corresponding singularities are constructed in Section \ref{sec.emb}
where we prove Theorem \ref{main.thm.1} except for an explicit resolution
of the resulting singularities which is accomplished in
Section \ref{sec.5}.

The proof of Theorem \ref{cm.main.thm.intro} is given in Section \ref{CM.sec}
where several other equivalent conditions are also treated.
 Theorem \ref{main.thm.2} on rational singularities
is reviewed in  Section \ref{Ratl.sec}.

Open questions and problems are  discussed in  Section \ref{sec.open}.

 \begin{ack}
I thank I.~Dolgachev, T.~de~Fernex, 
T.~Jarvis, M.~Kapovich, L.~Maxim, A.~N\'emethi, P.~Ozsv\'ath, S.~Payne, 
P.~Popescu-Pampu,
M.~Ramachand\-ran, 
J.~Shaneson, T.~Szamuely, D.~Toledo, J.~Wahl
and C.~Xu
 for comments and corrections.
Partial financial support   was provided  by  the NSF under grant number 
DMS-07-58275 and by the Simons Foundation. Part of this paper was written 
while the author visited the University of Utah.
\end{ack}

\section{Weighted homogeneous links}\label{sec.1}

\begin{defn}[Weighted homogeneous singularities]
Assign positive weights to the variables $w(x_i)\in \z$,
then the  weight of a monomial $\prod_i x_i^{a_i}$ is
$$
w\bigl(\tprod_i x_i^{a_i}\bigr):=\tsum_i a_i w(x_i).
$$
A polynomial $f$ is called {\it weighted homogeneous} of
weighted-degree $w(f)$ iff every monomial that occurs 
in $f$ with nonzero coefficient has weight $w(f)$. 

Fix weights ${\mathbf w}:=\bigl(w(x_1), \dots, w(x_N)\bigr)$ and let
$\{f_i:i\in I\}$ be weighted homogeneous polynomials.
They define both a projective variety in a weighted projective space
$$
Z(f_i:i\in I)\subset \p({\mathbf w})
$$
and an affine {\it weighted cone}
$$
C(f_i:i\in I)\subset \c^N.
$$ 

Somewhat loosely speaking, a singularity is called  {\it weighted homogeneous} 
if  it is isomorphic to a 
singularity defined by  a weighted cone 
for some weights $w(x_i)$.
(In the literature these are frequently called  {\it quasi-homogeneous} singularities.)

In many cases the weights are uniquely determined by the singularity
(up to rescaling) but not always. For instance, the
singularity  $(xy=z^n)$ is weighted homogeneous for any weights
that satisfy $w(x)+w(y)=n\cdot w(z)$.

If $C\subset \c^N$ is a weighted cone then
it has a $\c^*$-action given by
$$
(x_1,\dots, x_N)\mapsto \bigl(t^{m_1}x_1,\dots, t^{m_N}x_N\bigr)
\qtq{where} m_i=\tfrac1{w(x_i)}\tprod_j w(x_j).
$$
Conversely, let $X$ be a variety with a $\c^*$ action and
$x\in X$ a fixed point that is attractive as $t\to 0$.
Linearizing the action shows that $x\in X$ is a 
weighted homogeneous singularity.
\end{defn}

\begin{say}[Links of weighted homogeneous  singularities]
\label{link.of.qhs.say}
The $\c^*$-action on a weighted homogeneous singularity $(x\in X)$
induces a fixed point free $\s^1$-action on its link $L$.
If we think of $X$ as a weighted cone over the corresponding
projective variety  $Z\subset \p({\mathbf w})$
then we get a projection $\pi:L\to Z$ whose fibers are exactly the
orbits of the $\s^1$-action,
that is, the link of a weighted homogeneous singularity
has a {\it Seifert bundle} structure.
(For our purposes we can think that a Seifert bundle is the same as
a fixed point free $\s^1$-action.)
If  $(x\in X)$ is an isolated singularity then $Z$ is an orbifold.

It is thus natural to study the topology of 
links of  weighted homogeneous singularities in two steps.
\begin{enumerate}
\item Describe all $2n-1$-manifolds with a fixed point free $\s^1$-action.
\item Describe which among them occur as links of
weighted homogeneous singularities.
\end{enumerate}
\end{say}

\begin{say}[Homology of a weighted homogeneous  link]
\cite{or-wa}\label{seif.Q.homology.say}
Let $\pi:L\to Z$ be the Seifert bundle structure.
The cohomology of $L$ is computed by a spectral sequence
$$
H^i\bigl(Z, R^j\pi_*\q_L\bigr)\Rightarrow H^{i+j}(L, \q).
\eqno{(\ref{seif.Q.homology.say}.1)}
$$
All the fibers are oriented circles, thus
$R^0\pi_*\q_L\cong R^1\pi_*\q_L\cong\q_Z$ and $R^j\pi_*\q_L=0$ for $j>1$. 
Thus the $E_2$-term of the spectral sequence  is
$$
\begin{gathered}
    \xymatrix{%
H^0(Z, \q)\ar[drr] & H^1(Z, \q)\ar[drr]  & H^2(Z, \q) & \cdots\\
H^0(Z, \q) & H^1(Z, \q) & H^2(Z, \q) & \cdots
 }
  \end{gathered}
\eqno{(\ref{seif.Q.homology.say}.2)}
$$
where the differentials  are cup product with the (weighted) hyperplane class
$$
c_1\bigl(\o_Z(1)\bigr)\cup: H^i(Z, R^1\pi_*\q_L)\cong H^i(Z, \q)
\to H^{i+2}(Z, \q).
\eqno{(\ref{seif.Q.homology.say}.3)}
$$
Since $Z$ is an orbifold, these are injective if
$i+2\leq \dim Z$ and surjective if $i\geq \dim Z$.
Thus we conclude that
$$
\begin{array}{lcl}
h^i(L, \q)&=&h^i(Z, \q)- h^{i-2}(Z, \q)\qtq{if $i\leq \dim Z$ and}\\[1ex]
h^{i+1}(L, \q)&=&h^i(Z, \q)- h^{i+2}(Z, \q)\qtq{if $i\geq \dim Z$}
\end{array}
\eqno{(\ref{seif.Q.homology.say}.4)}
$$
where we set $h^i(Z, \q)=0$ for $i<0$ or $i>2\dim Z$.
In particular we see that $L$ is a rational homology sphere iff
$Z$ is a rational homology complex projective space.

By contrast, the spectral sequence computing the integral
cohomology of $L$ is much more complicated. We have a natural injection
$R^1\pi_*\z_L\into  \z_Z$ which is, however, rarely an isomorphism.
The computations were  carried out only for $\dim L\leq 5$
\cite{k-em2}. 
\end{say}

\begin{say}[Weighted homogeneous surface singularities]
This is the only case that is fully understood.

The classification of fixed point free circle actions on 3--manifolds
was considered by Seifert \cite{seif}. If $M$ is a  3--manifold
with a fixed point free circle action then the quotient space
$F:=M/\s^1$ is a surface (without boundary in the orientable case).
The classification of these {\it Seifert fibered} 
3--manifolds $f:M\to F$ is thus equivalent to the
 classification of fixed point free circle actions.
It should be noted that already in this classical case,
it is conceptually better to view the base surface $F$
not as a 2--manifold but as a 2-dimensional {\it orbifold},
see \cite{scott} for a detailed survey from this point of view.

Descriptions of weighted homogeneous surface singularities
are given in \cite{pinkham, dol-link, demazure, fl-za}.
\end{say}

\subsection*{Weighted homogeneous 3-fold singularities}{\ }

There is a quite clear picture about the simply connected case
since simply connected 5--manifolds are determined by their 
homology. 

By a theorem of  \cite{smale, barden}, 
 a  simply connected, compact 5--manifold $L$ is uniquely determined by 
$H_2(L,\z)$ and the second
Stiefel--Whitney class, which we view as a map
$w_2: H_2(L,\z)\to \z/2$.
Furthermore, there is such a 5--manifold 
iff there is an integer  $k\geq 0$ and a finite Abelian group $A$
such that either
 $H_2(L,\z)\cong \z^k+A+A$ and
$w_2: H_2(L,\z)\to \z/2$ is arbitrary, or 
 $H_2(L,\z)\cong \z^k+A+A+\z/2$ and
$w_2$ is projection on the $\z/2$-summand.

The existence of Seifert bundles on 
simply connected compact 5--manifolds was treated in \cite{k-circ}.
The answer mostly depends on the torsion subgroup of
 $H_2(L,\z)$, but there is a subtle interplay
with $w_2$.

\begin{defn}\label{i(L).defn}
Let $M$ be any manifold. Write its  second homology as
 a direct sum
of cyclic groups of prime power order
$$
H_2(M,\z)=\z^k+\tsum_{p,i} \bigl(\z/p^i\z\bigr)^{c(p^i)}
\eqno{(\ref{i(L).defn}.1)}
$$
for some $k=\dim H_2(M,\q)$ and $ c(p^i)=c(p^i,M)$.
The numbers $k, c(p^i)$ are  determined by
$H_2(M,\z)$ but the subgroups
$(\z/p^i)^{c(p^i)}\subset H_2(M,\z)$ are  usually not unique.
One can choose the decomposition
(\ref{i(L).defn}.1) such that 
$w_2:H_2(M,\z)\to \z/2$ is zero 
on all but one summand  $\z/2^n$. This value $n$ is unique
and it is   denoted by $i(M)$ \cite{barden}.
This invariant
 can take up any value $n$ for which $c(2^n)\neq 0$,
besides $0$ and $\infty$. 
Alternatively, 
$i(M)$ is the
smallest $n$ such that there is an $\alpha\in H_2(M,\z)$ such that
$w_2(\alpha)\neq 0$ and $\alpha$ has order $2^n$.
\end{defn}

The existence of a 
fixed point free differentiable circle action 
puts strong restrictions on $H_2$ and on $w_2$.

\begin{thm}\cite[Thm.3]{k-circ}\label{main.thm.2.k-circ}
Let  $L$ be a   compact, simply connected  5--manifold.
Then $L$ admits a fixed point free differentiable circle action
if and only if $H_2(L,\z)$ and $w_2$  satisfy
the following conditions.
\begin{enumerate}
\item  For every  $p$, we have
at most $\dim H_2(M,\q)+1$ nonzero $c(p^i)$ in 
(\ref{i(L).defn}.1). 
\item One can arrange that $w_2:H_2(L,\z)\to \z/2$
is the zero map on all but the $\z^k+(\z/2)^{c(2)}$ summands
in (\ref{i(L).defn}.1). That is, $i(L)\in \{0,1,\infty\}$.
\item If $i(L)=\infty$ then $\#\{i:c(2^i)>0\}\leq \dim H_2(M,\q)$.
\end{enumerate}
\end{thm}

\begin{rem} Note that while Theorem \ref{main.thm.2.k-circ} tells us which
compact, simply connected  5--manifolds admit a 
fixed point free differentiable circle action,
the proof does not classify all circle actions.
In particular, the classification of all circle actions
on $\s^5$ is not known.
\end{rem}

By contrast very little is known about which  
compact, simply connected  5--manifolds  occur as
links of 
weighted homogeneous singularities. It is known that not every
Seifert bundle occurs \cite[Lem.49]{k-circ} but
a full answer seems unlikely.

Nothing seems to be known in higher dimensions.

\begin{say}[Einstein metrics on weighted homogeneous links]

By a result of \cite{MR0154235}, the link of a cone over a 
smooth projective variety $Z\subset \p^N$ carries a natural Einstein metric
iff $-K_Z$ is a positive multiple of the hyperplane class and
$Z$ carries a K\"ahler--Einstein metric.
This was generalized by \cite{bg00} to  weighted cones.
Here one needs to work with an orbifold canonical class
$K_X+\Delta$ and a suitable orbifold K\"ahler--Einstein metric on $(X,\Delta)$.

This approach was used to construct new Einstein metrics on
spheres and exotic spheres \cite{bgk, bgkt} and on many
5-manifolds  \cite{k-em2, k-em1, k-em3}.

See \cite{bg-book} for a comprehensive treatment.

\end{say}

\section{Construction of singularities}\label{sec.2}

The  construction has 5 main steps, none of which 
is fully understood at the moment.
After summarizing them, we discuss the difficulties in  more
detail.
Although the steps can not be carried out
in full generality, we understand enough about them to obtain
the main theorems.

\begin{say}[Main steps of the construction]\label{main.steps.say}{\ }

{\it Step.}\ref{main.steps.say}.1. For a simplicial complex  $C$
construct projective simple normal crossing varieties $V(C)$
such that ${\mathcal D}\bigl(V(C)\bigr)\cong C$.
\medskip

{\it Step.}\ref{main.steps.say}.2. 
 For a projective simple normal crossing variety $V$ construct
a smooth variety $Y(V)$ that contains $V$ as a divisor.
\medskip

{\it Step.}\ref{main.steps.say}.3. For a smooth variety $Y$
containing a simple normal crossing divisor $D$
construct an isolated singularity $(x\in X)$
such that $(D\subset Y)$ is a resolution of $(x\in X)$.
\medskip

{\it Step.}\ref{main.steps.say}.4. Describe the link
$L(x\in X)$ in terms of the topology of $D$ and the
Chern class of the normal bundle of $D$.
\medskip

{\it Step.}\ref{main.steps.say}.5. Describe the
relationship between the properties of the singularity
 $(x\in X)$ and the original  simplicial complex  $C$.

\end{say}

\begin{say}[Discussion of  Step \ref{main.steps.say}.1]
\label{disc.main.steps.say.1}
 I believe that for every simplicial complex  $C$
there are many projective simple normal crossing varieties $V(C)$
such that ${\mathcal D}\bigl(V(C)\bigr)\cong C$.\footnote{This is now proved in 
\cite{k-dual2}.} 

There seem to be two main difficulties of a step-by-step approach.

First, topology would suggest that one should build up the skeleta
of $V(\CC)$  one dimension at a time. It is easy to obtain the 1-skeleton
by gluing rational curves. The 2-skeleton is still straightforward
since rational surfaces do contain cycles of rational curves of 
arbitrary length. However, at the next step we run into a problem
similar to  Step \ref{main.steps.say}.2 and usually a 2-skeleton
can not be extended to a 3-skeleton.  Our solution in \cite{k-fg}
is to work with triangulations of $n$-dimensional submanifolds with
boundary in $\r^n$. The ambient $\r^n$ gives a rigidification
and this makes it possible to have a consistent choice for all the strata.

Second, even if we construct a simple normal crossing variety $V$,
it is not easy to decide whether it is projective.
This is illustrated by the following example of 
``triangular pillows'' \cite[Exmp.34]{k-fg}.

Let us start with an example that is not simple normal crossing. 

Take 2 copies  $\p^2_i:= \p^2(x_i:y_i: z_i)$ of $\c\p^2$ and 
the triangles  $C_i:=(x_iy_iz_i=0)\subset \p^2_i$.
Given $c_x, c_y, c_z\in \c^*$, 
define $\phi(c_x, c_y, c_z):C_1\to C_2$ by 
$(0:y_1:z_1)\mapsto (0:y_1:c_z z_1)$, $(x_1:0:z_1)\mapsto (c_xx_1:0: z_1)$
and $(x_1:y_1:0)\mapsto (x_1:c_yy_1: 0)$ 
and glue the 2 copies  of $\p^2$ 
using  $\phi(c_x, c_y, c_z)$ to get the surface
$S(c_x, c_y, c_z)$. 

We claim that $S(c_x, c_y, c_z)$
 is projective iff the product $c_xc_yc_z$ is a root of unity.

To see this  note that $\pic^0(C_i)\cong \c^*$
and  $\pic^r(C_i)$ is a principal homogeneous space
under $\c^*$ for every $r\in \z$. We can identify $\pic^3(C_i)$
 with $\c^*$ using
the restriction of the ample generator $L_i$ of 
$\pic\bigl(\p^2_i\bigr)\cong \z$ as the base point.

The key observation is that $\phi(c_x, c_y, c_z)^*:\pic^3(C_2)\to 
\pic^3(C_1)$ is 
multiplication by $c_xc_yc_z$. 
Thus if $c_xc_yc_z$ is an $r$th root of unity then 
$L_1^r$ and $L_2^r$ glue together to an ample line bundle
 but otherwise $S(c_x, c_y, c_z)$ carries only the  trivial line bundle.

We can create a similar 
simple normal crossing example by smoothing the triangles $C_i$.
That is, we take 2 copies  $\p^2_i:= \p^2(x_i:y_i: z_i)$ of $\c\p^2$ and 
smooth elliptic curves   $E_i:=(x_i^3+y_i^3+z_i^3=0)\subset \p^2_i$.

Every automorphism $\tau\in \aut(x^3+y^3+z^3=0)$
can be identified with an isomorphism $\tau:E_1\cong E_2$,
giving a 
simple normal crossing surface  $S(\tau)$. The above argument then shows that
$S(\tau)$ is projective iff $\tau^m=1$ for some $m>0$.

These examples are actually not surprising. One can think of the
surfaces $S(c_x, c_y, c_z)$ and $S(\tau)$ as  degenerate K3
surfaces of degree 2 and K3 surfaces have non-projective deformations.
Similarly, $S(c_x, c_y, c_z)$ and $S(\tau)$ can be non-projective.
One somewhat unusual aspect is that while a smooth K3 surface
is projective iff it is a scheme,  the above singular examples
are always schemes yet many of them are  non-projective.

\end{say}

\begin{say}[Discussion of  Step \ref{main.steps.say}.2]
\label{disc.main.steps.say.2}
This is surprisingly subtle. First note that not every
projective simple normal crossing variety $V$ can be realized as a
divisor on 
a smooth variety $Y$. A simple obstruction is the following.

Let $Y$ be a  smooth variety and $D_1+D_2$ a 
simple normal crossing divisor on $Y$.
Set $Z:=D_1\cap D_2$. Then $N_{Z,D_2}\cong N_{D_1, Y}|_Z$
where $N_{X,Y}$ denotes the normal bundle of $X\subset Y$.

Thus if $V=V_1\cup V_2$ is a 
simple normal crossing variety with $W:=V_1\cap V_2$
such that $N_{W,V_2}$ is not the restriction of any line bundle
from $V_1$ then $V$ is not  a
simple normal crossing divisor in a  nonsingular variety.

I originally hoped that such normal bundle considerations give
necessary and sufficient conditions, but
recent examples of \cite{2012arXiv1206.1994F, 2012arXiv1206.2475F}
 show that this is not the case.

For now, no necessary and sufficient conditions of embeddability are known.
In the original papers 
\cite{k-t-lc-exmp, k-fg}
we went around this problem
by first embedding a simple normal crossing variety $V$
into a singular variety $Y$ and then showing that for
the purposes of computing the  fundamental group of the link
the singularities of $Y$ do not matter.

We improve on this in Section \ref{sec.5}.
\end{say}

\begin{say}[Discussion of  Step \ref{main.steps.say}.3]
\label{disc.main.steps.say.3}
By a result of \cite{artin}, a compact divisor 
contained in a smooth variety $D=\cup_i D_i\subset Y$ can be contracted to
a point if there are positive integers $m_i$ such that
$\o_Y(-\sum_i m_i D_i)|_{D_j}$ is ample for every $j$. 

It is known that this condition is not necessary and no 
necessary and sufficient characterizations are known. 
However, it is easy to check the above condition 
in our examples.
\end{say}

\begin{say}[Discussion of  Step \ref{main.steps.say}.4]
\label{disc.main.steps.say.4}
This approach, initiated in \cite{mumf-top}, 
 has been especially successful for surfaces.

In principle the method of  \cite{mumf-top} leads to a
complete description of the link, but it seems rather
difficult to perform explicit computations. Computing the
fundamental group of the links seems rather daunting in general.
Fortunately, we managed to find some simple conditions
that ensure that the natural maps
$$
\pi_1\bigl(L(x\in X)\bigr)\to \pi_1\bigl(\res(X)\bigr)\to 
\pi_1\bigl(\dres(X)\bigr)
$$
are isomorphisms. 
However, these  simple conditions force $D$ to be more complicated
than necessary, in particular we seem to lose control of the
canonical class of $X$.
\end{say}

\begin{say}[Discussion of  Step \ref{main.steps.say}.5]
\label{disc.main.steps.say.5} 

For surfaces there is a very tight connection
between the topology of the link and the
algebro-geometric properties of a singularity.
In higher dimension, one can obtain very little information
from the topology alone. As we noted, there are many examples
where $X$ is a topological manifold yet very singular
as a variety.

There is more reason to believe that algebro-geometric properties
restrict the topology. For example, the results of
Section \ref{CM.sec} rely on the observation that if
$(x\in X)$ is a rational (or even just 1-rational) singularity
then $H_1\bigl(L(x\in X), \q\bigr)=0$.
\end{say}

\section{Voronoi complexes}
\label{sec.vor}

\begin{defn}
A (convex) {\it Euclidean polyhedron} is a subset $P$ of $\r^n$ given
 by a finite collection of linear inequalities 
(some of which may be strict and some not). 
A {\em face} of $P$ is a subset of $P$ which is given by converting some of 
these non-strict inequalities to equalities.

A {\em Euclidean polyhedral complex} in $\r^n$  is a 
collection of closed Euclidean polyhedra $\CC$ in  $\r^n$
such that
\begin{enumerate}
\item if $P\in \CC$  then every face of $P$ is in $\CC$ and
\item if $P_1, P_2\in \CC$ then $P_1\cap P_2$
is a face of both of the $P_i$ (or empty).
\end{enumerate}
The union of the faces of a Euclidean polyhedral complex
$\CC$ is denoted by $|\CC|$. 
\end{defn}

For us the most important examples are the following.

\begin{defn}[Voronoi complex]\label{voronoi.defn}
Let $Y=\{y_i:i\in I\}\subset \r^n$ be a finite subset. 
For each $i\in I$  the corresponding  {\it Voronoi cell} is 
the set of points that are closer to $y_i$ than to any other $y_j$,
that is
$$
V_i:=\{x\in \r^n: d(x,y_i)\le d(x,y_j), \forall j\in I\}
$$
where $d(x,y)$ denotes the Euclidean distance.
Each cell $V_i$ is a closed (possibly unbounded) polyhedron in $\r^n$.

The Voronoi cells and their faces give a 
Euclidean polyhedral complex, called the {\it Voronoi complex}
or {\it Voronoi tessellation}
associated to $Y$.

For a subset $J\subset I$ let $H_J$ denote the linear subspace
$$
H_J:=\{ x\in \r^n: d(x,y_i)= d(x,y_j)\ \forall i,j\in J\}.
$$
The affine span of each face of the Voronoi complex is 
one of the $H_J$. If $J$ has 2 elements $\{i,j\}$
then $H_{ij} $ is a hyperplane
$H_{ij}=\{ x\in \r^n: d(x,y_i)= d(x,y_j)\} $.

A  Voronoi complex  is called {\it simple} if for every $k$, 
every codimension $k$ face is contained in exactly $k+1$
Voronoi cells. 
Not every Voronoi complex is simple, but it is easy to see that
among  finite subsets  $Y\subset \r^n$ those with a 
simple Voronoi complex $\CC(Y)$ form an 
 open and dense set.

Let $\CC$ be a  simple Voronoi complex.
For each face $F\in \CC$, let
$V_i$ for $i\in I_F$ be the Voronoi cells containing $F$.
The vertices $ \{y_i:i\in I_F\}$ form a simplex whose dimension
equals the codimension of $F$. These simplices define the
 {\it Delaunay  triangulation} dual to $\CC$.
\end{defn}

\begin{thm}\cite[Cor.21]{k-fg}\label{cor.21.thm}
 Let $T$ be a finite simplicial complex of dimension $n$.
Then there is an embedding $j:T\into \r^{2n+1}$,
a simple Voronoi complex $\CC$ in $ \r^{2n+1}$ and
a subcomplex
 $\CC(T)\subset \CC$   of pure dimension $2n+1$ containing $j(T)$ such that
 the inclusion
$j(T)\subset |\CC(T)|$ is a  homotopy equivalence.
\end{thm}

Outline of the proof. First we embed $T$ into $\r^{2n+1}$.
This is where the dimension increase comes from. 
 (We do not need an actual embedding, only
an embedding up-to homotopy, which is usually easier to get.)

Then we first use a result of \cite{Hirsch}
which says that if  $T$ is a finite simplicial complex in a smooth manifold 
${\mathbf R}$ then there exists a codimension 0 compact submanifold 
$M\subset {\mathbf R}$ with 
smooth boundary containing $T$ such that the inclusion
$T\subset M$ is a  homotopy equivalence.

Finally we construct a Voronoi complex using $M$.

Let $M\subset \r^m$ be a compact subset, 
 $Y\subset \r^m$  a finite set of points
and $\CC(Y)$ the corresponding Voronoi complex.
Let $\CC_m(Y,M)$ be the collection of those $m$-cells 
in the Voronoi complex $\CC(Y)$
whose intersection
with $M$ is not empty  and $\CC(Y,M)$ the polyhedral complex consisting of
the cells in  $\CC_m(Y,M)$ and their faces. 
Then $M\subset |\CC(Y,M)|$.

We conclude by using a theorem of \cite{Cairns} that says that
if $M$ is a $C^2$-submanifold with $C^2$-boundary
then for a 
 suitably fine mesh of points $Y\subset \r^{m}$
the inclusion $M\subset |\CC(Y,M)|$ is a homotopy equivalence. \qed

\section{Simple normal crossing varieties}
\label{sec.snc}

Let $\CC$ be a 
purely $m$-dimensional, compact  subcomplex of a 
simple Voronoi complex  in $\r^m$.
Our aim is to construct a projective simple normal crossing variety
$V(\CC)$ whose dual complex naturally identifies with
 the Delaunay  triangulation of $\CC$.

\begin{say}[First attempt]\label{first.att.say} For each $m$-polytope
$P_i\in \CC$ we associate a copy   $\p^m_{(i)}=\c\p^m$.
For a subvariety $W\subset \c\p^m$ we let
$W_{(i)}$ or $W^{(i)}$ denote the corresponding subvariety of $\p^m_{(i)}$.

If $P_i$ and $P_j$ have a common face $F_{ij}$ of dimension $m-1$
then
the complexification of the affine span of $F_{ij}$ gives
hyperplanes  $H^{(i)}_{ij}\subset \p^m_{(i)}$ and 
$H^{(j)}_{ij}\subset \p^m_{(j)}$. Moreover, $H^{(i)}_{ij}$ and $H^{(j)}_{ij}$
come with a natural identification
$\sigma_{ij}:H^{(i)}_{ij}\cong H^{(j)}_{ij}$.

We use $\sigma_{ij}$ to glue  $\p^m_{(i)}$ and  $\p^m_{(j)}$ together.
The resulting variety is isomorphic to the union of 2 hyperplanes
in $\c\p^{m+1}$.

It is harder to see what happens if we try to perform all these
gluings $\sigma_{ij}$ simultaneously.

Let $\amalg_i \p^m_{(i)}$ denote the disjoint union of all
the $\p^m_{(i)}$. Each $\sigma_{ij}$ defines a relation that
identifies a point  $p_{(i)}\in H^{(i)}_{ij}\subset \p^m_{(i)}$
with its image  $p_{(j)}=\sigma_{ij}(p_{(i)})\in H^{(j)}_{ij}\subset \p^m_{(j)}$.
Let $\Sigma$ denote the equivalence relation generated by all the
$\sigma_{ij}$. 

It is easy to see 
(cf.\ \cite[Lem.17]{k-q})
 that there is a projective algebraic variety
$$
 \amalg_i \p^m_{(i)} \longrightarrow 
\bigl(\amalg_i \p^m_{(i)}\bigr)/\Sigma \longrightarrow 
\c\p^m
$$
whose points are exactly the equivalence classes of $\Sigma$.

This gives the correct simple normal crossing variety if $m=1$ but already
for $m=2$ we have problems. 
For instance, consider three 2-cells
$P_i, P_j, P_k$ such that $P_i$ and $ P_j$ have a common face $F_{ij}$,
$P_j$ and $ P_k$ have a common face $F_{jk}$ but 
$P_i\cap  P_k=\emptyset$.
The problem is that while $F_{ij}$ and $F_{jk}$ are disjoint,
their complexified spans are lines in $\c\p^2$ hence they
intersect at a point $q$. Thus
$\sigma_{ij}$ identifies $q_{(i)}\in \p^2_{(i)}$ with $q_{(j)}\in \p^2_{(j)}$
and $\sigma_{jk}$ identifies $q_{(j)}\in \p^2_{(j)}$ with $q_{(k)}\in \p^2_{(k)}$
thus the equivalence relation $\Sigma$ identifies
 $q_{(i)}\in \p^2_{(i)}$ with $q_{(k)}\in \p^2_{(k)}$.
Thus in $\bigl(\amalg_i \p^m_{(i)}\bigr)/\Sigma $
the images of $\p^2_{(i)}$ and of $\p^2_{(k)}$ are not disjoint.

In order to get the correct simple normal crossing variety,
we need to remove these extra intersection points.
In higher dimensions we need to remove various linear subspaces as well. 
\end{say}

\begin{defn}[Essential and parasitic intersections] 
Let $\CC$ be a Voronoi complex on $\r^m$
defined by the points $\{y_i:i\in I\}$.
We have the linear subspaces $H_J$ defined in (\ref{voronoi.defn}). 
Assume for simplicity that $J_1\neq J_2$ implies that
$H_{J_1}\neq H_{J_2}$.  

Let $P\subset \r^m$ be a  Voronoi cell.
We say that $H_J$ is {\it essential} for $P$ 
if it is the affine span of a face of $P$.
Otherwise it is called {\it parasitic} for $P$. 
\end{defn}

\begin{lem} Let $P\subset \r^m$ be a simple Voronoi cell.
\begin{enumerate}
\item Every essential subspace $L$ of dimension $\leq m-2$ is contained
in a unique smallest  parasitic subspace which has dimension $\dim L+1$.
\item The intersection of two  parasitic subspaces is again parasitic.
\end{enumerate}
\end{lem}

Proof. There is a  point $y_p\in P$ and a subset $J\subset I$
such that  $H_{ip}$ are spans of faces of $P$ for $i\in J$ and
$L=\cap_{i\in J}  H_{ip}$. Thus the unique
$\dim L+1$-dimensional parasitic subspace containing $L$ is
$H_J$.

Assume that $L_1,L_2$ are parasitic.
If $L_1\cap L_2$ is essential then there is a unique
smallest parasitic subspace $L'\supset L_1\cap L_2$. Then
$L'\subset L_i$ a contradiction. \qed

\begin{say}[Removing parasitic intersections]\label{rem.par.say}
Let $\{H_s:s\in S\}$ be a finite set of hyperplanes of $\c\p^m$.
For $Q\subset S$ set $H_Q:=\cap_{s\in Q} H_s$. 
Let ${\mathcal P}\subset 2^S$ be a subset closed under unions. 

 Set  $\pi_0:P^0\cong \c\p^m$.  
If $\pi_r:P^r\to \c\p^m$ is already defined then let
$P^{r+1}\to P^r$ denote the blow-up of the union of
birational transforms of all the $H_Q$ such that $Q\in {\mathcal P}$
and $\dim H_Q=r$. Then $\pi_{r+1}$ is the composite
$P^{r+1}\to P^r\to \c\p^m$.

Note that we blow up a disjoint union of
smooth subvarieties since any intersection of the
$r$-dimensional $H_Q$ is lower dimensional, hence
it was removed by an earlier blow up.
Finally set  $\Pi:\tilde P:=P^{m-2}\to \c\p^m$.
\end{say}

Let $\CC$ be a pure dimensional subcomplex of a Voronoi complex
as in (\ref{cor.21.thm}). For each cell $P_i\in \CC$
we use (\ref{rem.par.say}) with
$$
{\mathcal P}_i:=\{\mbox{parasitic intersections for $P_i$}\}
$$
to obtain $\tilde P_{(i)} $. Note that if $P_i$ and $P_j$ have a common 
codimension 1 face $F_{ij}$ then we perform the same blow-ups on
the complexifications $H_{ij}^{(i)}\subset \p^m_{(i)}$ and 
$H_{ij}^{(j)}\subset \p^m_{(j)}$. 
 Thus $\sigma_{ij}:H_{ij}^{(i)}\cong H_{ij}^{(j)}$ lifts to
the birational transforms
$$
\tilde \sigma_{ij}:\tilde H_{ij}^{(i)}\cong \tilde H_{ij}^{(j)}.
$$
As before, the $\tilde \sigma_{ij} $ define an
equivalence relation $\tilde \Sigma$ on $ \amalg_i \tilde P_{(i)}$.
 With these changes,
the approach outlined in Paragraph \ref{first.att.say}
 does work and we get the following.

\begin{thm}\cite[Prop.28]{k-fg} \label{snc.exists.main.thm}
With the above notation there is a
projective, simple normal crossing  variety 
$$
V(\CC):=\bigl(\amalg_i \tilde P_{(i)}\bigr)/\tilde \Sigma
$$
with the following properties.
\begin{enumerate}
\item There is a  finite morphism
 $\amalg_i \tilde P_{(i)} \longrightarrow V(\CC)$
whose fibers are exactly the equivalence classes of $\tilde \Sigma$.
\item  The dual complex ${\mathcal D}\bigl(V(\CC)\bigr)$ is naturally
identified with the  Delaunay  triangulation of 
$\CC$.
\end{enumerate}
\end{thm}

Comments on the proof. The existence of $V(\CC)$ is relatively easy
either directly as in \cite[Prop.31]{k-fg} or using the general theory
of quotients by finite equivalence relations as in \cite{k-q}.

As we noted in Paragraph \ref{disc.main.steps.say.1}
the projectivity of such quotients is a rather delicate question
since the maps $ \tilde P_{(i)}\to \c\p^m$ are not finite any more.

The main advantage we have here is that each $\tilde P_{(i)}$
comes with a specific sequence of blow-ups
$\Pi_i:\tilde P_{(i)}\to \c\p^m$ and this enables us to write down
explicit, invertible, ample  subsheaves
$A_i\subset \Pi_i^*\o_{\c\p^m}(N)$ for some $N\gg 1$
that glue together to give an ample invertible
 sheaf on $V(T)$. For details see
\cite[Par.32]{k-fg}.\qed 
\medskip

The culmination of the results of the last 2 sections is
the following.

\begin{thm}\cite[Thm.29]{k-fg}\label{k-fg.thm}
 Let $T$ be a finite cell complex. Then there is a projective
simple normal crossing variety $Z_T$ such that
\begin{enumerate}
\item ${\mathcal D}(Z_T)$ is homotopy equivalent to $T$,
\item $\pi_1(Z_T)\cong \pi_1(T)$ and
\item $H^i\bigl(Z_T, \o_{Z_T}\bigr)\cong H^i(T, \c)$ for every $i\geq 0$.
\end{enumerate}
\end{thm}

Proof. We have already established (1) in (\ref{snc.exists.main.thm}),
 moreover the construction yields a  simple normal crossing variety $Z_T$
whose strata are all rational varieties.
In particular  every stratum $W\subset Z_T$ is simply connected and 
$H^r\bigl(W, \o_{W}\bigr)=0$  for every $r>0$.
Thus (2--3) follow from Lemmas \ref{friedman.lem}--\ref{pi_1-cong}. \qed
\medskip

The proof of the following lemma is essentially in 
\cite[pp.68--72]{gri-sch}. More explicit versions can be found in
\cite[pp.26--27]{friedman-etal} and 
\cite{Ishii85, 2009arXiv0902.4234A}.

\begin{lem}\label{friedman.lem}
Let $X$ be a simple normal crossing variety over $\c$ with
irreducible components $\{X_i:i\in I\}$. Let $T=D(X)$ be the dual 
 complex of $X$. 
\begin{enumerate}
\item There are natural injections
$H^r\bigl(T, \c\bigr)\into H^r\bigl(X, \o_{X}\bigr)$ for every $r$.
\item 
Assume that   $H^r\bigl(W, \o_{W}\bigr)=0$  for every $r>0$ and
for every stratum $W\subset X$.
Then  $H^r\bigl(X, \o_{X}\bigr)=H^r\bigl(T, \c\bigr)$ for every $r$.\qed
\end{enumerate}
\end{lem}

The following comparison result is rather straightforward.

\begin{lem}\label{pi_1-cong}\cite[Prop.3.1]{Corson}
Using the notation of (\ref{friedman.lem}) assume that
every stratum $W\subset X$ is $1$-connected. Then 
 $\pi_1(X)\cong \pi_1\bigl({\mathcal D}(X)\bigr)$. \qed
\end{lem}

\section{Generic embeddings of 
simple normal crossing varieties}
\label{sec.emb}

The following is a summary of the construction of \cite{k-t-lc-exmp};
see also \cite[Sec.3.4]{kk-singbook} for an improved version.

\begin{say}\label{summary.k-exmp}
Let $Z$ be a projective,  local complete intersection variety
of dimension $n$ 
and choose any embedding  $Z\subset P$ into a smooth projective 
variety of dimension $N$.  
(We can take  $P=\p^N$ for $N\gg 1$.)
 Let $L$ be a sufficiently ample line bundle on $P$.
Let $Z\subset Y_1\subset P$ be the  complete intersection
of $(N-n-1)$  general sections of $L(-Z)$. Set
$$
Y:= B_{(-Z)}Y_1:=\proj_{Y_1}\tsum_{m=0}^{\infty}\o_{Y_1}(mZ).
$$
(Note that this is not the
 blow-up of $Z$ but  the
 blow-up of its inverse in the  class group.)

It is proved in \cite{k-t-lc-exmp} that 
the birational transform of $Z$ in $Y$ is  a Cartier divisor isomorphic to $Z$
and there is a  contraction morphism
$$
\begin{array}{ccc}
Z & \subset & Y\\
\downarrow && \hphantom{\pi}\downarrow\pi \\
0 & \in & X
\end{array}
\eqno{(\ref{summary.k-exmp}.1)}
$$
such that $Y\setminus Z\cong X\setminus \{0\}$.
If $Y$ is smooth then  $\dres(0\in X)={\mathcal D}(Z)$ and we are done
with  Theorem \ref{main.thm.1}.
However, the construction of \cite{k-t-lc-exmp} yields a smooth variety
$Y$  only if $\dim Z=1$ or $Z$ is smooth.
(By (\ref{disc.main.steps.say.2}) this limitation is
not unexpected.)

In order to  resolve  singularities of $Y$
we need a  detailed description of them.
This is a local question, so we may assume that
 $Z\subset \c^N_{\mathbf x}$ is a complete intersection
  defined
by $f_1=\dots=f_{N-n}=0$. Let $Z\subset Y_1\subset \c^N$ be a general
complete intersection defined
by   equations
$$
h_{i,1} f_1+  \cdots  +h_{i,N-n} f_{N-n}  =  0
\qtq{for $i=1,\dots, N-n-1$. }
$$
Let $H=(h_{ij})$ be the $(N-n-1)\times (N-n)$ matrix of the system and 
$H_i$ the submatrix obtained by removing the $i$th column.
By \cite{k-t-lc-exmp} or \cite[Sec.3.2]{kk-singbook},
an open neighborhood of $Z\subset Y$ 
is defined by the equations
$$
\bigl(f_i=(-1)^i\cdot t\cdot \det H_i : i=1,\dots, N-n\bigr)
\subset \c^N_{\mathbf x}\times \c_t.
\eqno{(\ref{summary.k-exmp}.2)}
$$
Assume now that $Z$ has hypersurface singularities.
Up-to permuting the  $f_i$ and passing to a smaller open set, 
we may assume that
 $df_2,\dots, df_{N-n}$ are linearly independent everywhere along $Z$.
Then the singularities of $Y$ all come from the equation
$$
f_1=- t\cdot \det H_1.  
\eqno{(\ref{summary.k-exmp}.3)}
$$
Our aim is to write down  local normal forms for $Y$ along $Z$
in the normal crossing case.

On $\c^N$ there is a stratification
$\c^N=R_0\supset R_1\supset \cdots$
where $R_i$ is the set of points where
$\rank H_1\leq (N-n-1)-i$. Since the $h_{ij}$ are general,
$\codim_W R_i=i^2$ and we may assume that every stratum of $Z$
is transversal to each $R_i\setminus R_{i+1}$ (cf.\ Paragraph \ref{detvar.say}).

Let $S\subset Z$ be any stratum and $p\in S$ a point
such that $p\in R_m\setminus R_{m+1}$.
We can choose local coordinates  $\{x_1,\dots, x_d\}$ and 
$\{y_{rs}: 1\leq r,s\leq m\}$  such that, in a neighborhood of $p$,
$$
f_1=x_1\cdots x_d\qtq{and} 
\det H_1=\det\bigl( y_{rs}: 1\leq r,s\leq m\bigr).
$$
Note that $m^2\leq \dim S=n-d$, thus we can add $n-d-m^2$ further
coordinates $y_{ij}$ to get a complete local coordinate system on $S$.

Then  the $n$ coordinates $\{x_k,y_{ij}\}$ determine a 
map
$$
\sigma:\c^N\times \c_t\to \c^n  \times \c_t
$$
such that $\sigma(Y)$ is defined by the equation
$$
x_1\cdots x_d=t\cdot \det\bigl( y_{rs}: 1\leq r,s\leq m\bigr).
$$
Since  $df_2,\dots, df_{N-n}$ are linearly independent along $Z$,
we see that   $\sigma|_Y$ is \'etale along $Z\subset Y$.
\end{say}

We can summarize these considerations as follows.

\begin{prop}\label{k-esmp.prop}
 Let $Z$ be a normal crossing variety of dimension $n$. 
Then there is
a normal singularity  $(0\in X)$ of dimension $n+1$ and a proper, birational
morphism $\pi:Y\to X$ such that $\red\pi^{-1}(0)\cong Z$
and  for every point $p\in \pi^{-1}(0)$
we can choose local  (\'etale or analytic) coordinates
 called 
$\{x_i: i\in I_p\}$ and $ \{y_{rs}:1\leq r, s\leq m_p\}$
(plus possibly other unnamed coordinates)
such that one can write the local equations of $Z\subset Y$ as 
$$
\bigl( \tprod_{i\in I_p}x_i=t=0\bigr)\subset \Bigl(\tprod_{i\in I_p}x_i=
t\cdot \det\bigl(y_{rs}:1\leq r, s\leq m_p\bigr)\Bigr)
\subset \c^{n+2}.\qed
$$
\end{prop}

\begin{say}[Proof of Theorem \ref{main.thm.1}] \label{Ri=Hi.say} 
Let $T$ be a finite cell complex. By (\ref{k-fg.thm}) 
 there is a projective
simple normal crossing variety $Z$ such that
 ${\mathcal D}(Z)$ is homotopy equivalent to $T$,
 $\pi_1(Z)\cong \pi_1(T)$ and
 $H^i(Z, \o_{Z})\cong H^i(T, \c)$ for every $i\geq 0$.

Then Proposition \ref{k-esmp.prop} constructs a singularity
$(0\in X)$ with a partial resolution
$$
\begin{array}{ccc}
Z & \subset & Y\\
\downarrow && \hphantom{\pi}\downarrow\pi \\
0 & \in & X
\end{array}
\eqno{(\ref{Ri=Hi.say}.1)}
$$
The hardest is to check that we can resolve the singularities of
$Y$ without changing the homotopy type of the dual complex of
the exceptional divisor. This is done in Section \ref{sec.5}.

In order to show (\ref{main.thm.1}.2--3) we need
 further information about the varieties and
maps in (\ref{Ri=Hi.say}.1).

First, $Y$ has rational singularities. This is easy
to read off from their equations.
(For the purposes of Theorem \ref{link.thm}, we only need the case
$\dim Y=3$ when the only singularities we have are
 ordinary double points with local equation $x_1x_2=ty_{11}$.)

Second, we can arrange that $Z$ has very negative normal bundle in
$Y$. By a general argument this implies that
$R^i\pi_*\o_Y\cong H^i(Z, \o_Z)$, proving
(\ref{main.thm.1}.3);
see \cite[Prop.9]{k-t-lc-exmp} for details.

Finally we need to compare $\pi_1(Z)$ with
$\pi_1\bigl(L(0\in X)\bigr)$.  There is always a surjection
$$
\pi_1\bigl(L(0\in X)\bigr)\onto \pi_1(Z)
\eqno{(\ref{Ri=Hi.say}.2)}
$$
but it can have a large kernel. We claim however, that with suitable 
choices we can arrange that (\ref{Ri=Hi.say}.2) is an isomorphism.
It is easiest to work not on $Z\subset Y$ but on a resolution
$Z'\subset Y'$.

More generally, 
let $W$ be a smooth variety,
$D=\cup_iD_i\subset W$ a simple normal crossing divisor and
 $T\supset D$  a regular neighborhood with boundary
$M=\partial T$.
There is  a natural (up to homotopy) retraction map
$T\to D$ which induces $M\to D$ hence a surjection
$\pi_1(M)\onto \pi_1(D)$ whose kernel is generated
(as a normal subgroup) by  the simple loops $\gamma_i$ around the $D_i$.

In order to understand this kernel, assume first that
$D$ is smooth. Then $M\to D$ is a circle bundle hence there is
an exact sequence
$$
\pi_2(D)\stackrel{c_1\cap}{\longrightarrow}
\z\cong \pi_1(\s^1)\to \pi_1(M)\to \pi_1(D)\to 1
$$
where $c_1$ is the Chern class of the normal bundle
of $D$ in $X$. 
Thus if $c_1\cap\alpha=1$ for some $\alpha\in \pi_2(D)$
then $\pi_1(M)\cong \pi_1(D)$.
In the general case, arguing as above we see that
$\pi_1(M)\cong \pi_1(D)$ if the following holds:
\begin{enumerate}\setcounter{enumi}{2}
\item For every $i$ there is a class
$\alpha_i\in \pi_2\bigl(D_i^0\bigr)$ such that
$c_1\bigl(N_{D_i,X}\bigr)\cap\alpha_i=1$
where $D_i^0:=D_i\setminus\{\mbox{other components of $D$}\}$.
\end{enumerate}

Condition (3) is typically very easy to achieve in our constructions.
Indeed, we obtain the $D_i^0$ by starting with $\c\p^m$,
blowing it up many times and then removing  a few divisors.
Thus we end up with very large $H_2\bigl(D_i^0,\z\bigr)$
and typically the $D_i^0$ are even simply connected,
hence $\pi_2\bigl(D_i^0)=H_2\bigl(D_i^0,\z\bigr)$.\qed
\end{say}

\begin{say}[Determinantal varieties]\label{detvar.say}
We have used the following basic properties of determinantal varieties.
These are quite easy to prove directly; see \cite[12.2 and 14.16]{Harris95}
for a more general case.

Let $V$ be a smooth, affine variety,
and ${\mathcal L}\subset \o_V$
a finite dimensional sub vector space without common zeros.
Let  $H=\bigl(h_{ij}\bigr)$ be an $n\times n$ matrix
whose entries are general elements in ${\mathcal L}$.
For a point $p\in V$ set $m_p=\operatorname{corank} H(p)$. 
 Then there are local analytic coordinates
 $ \{y_{rs}:1\leq r, s\leq m_p\}$
(plus possibly other unnamed coordinates)
such that, in a neighborhood of $p$, 
$$
\det H= \det\bigl(y_{rs}:1\leq r, s\leq m_p\bigr).
$$
In particular, $\mult_p(\det H)=\operatorname{corank} H(p)$,
for every $m$  the set of points $R_m\subset V$ where
$\operatorname{corank} H(p)\geq m$ is a  subvariety
of pure codimension $m^2$ and $\sing R_m=R_{m+1}$.
\end{say}

\section{Resolution of generic embeddings}
\label{sec.5}

In this section we start with the varieties constructed in  
Proposition \ref{k-esmp.prop} and
 resolve their singularities.
Surprisingly, the resolution process 
described in Paragraphs \ref{indestup}--\ref{binres.say} 
 leaves the  dual 
complex unchanged and we get the following.

\begin{thm}\label{main.thm.3} Let $Z$ be a projective  simple normal crossing
 variety of dimension $n$.
Then there is a normal singularity $(0\in X)$ of dimension $(n+1)$
and a resolution  $\pi:Y\to X$ such that
$E:=\pi^{-1}(0)\subset Y$ is a  simple normal crossing divisor 
and  its dual  complex  ${\mathcal D}(E)$
is naturally identified with   ${\mathcal D}(Z)$.
(More precisely, there is a morphism $E\to Z$
that induces a birational map on every stratum.)
\end{thm}

\begin{say}[Inductive set-up for resolution]\label{indestup}
 The  object we try to resolve is  a triple
$$
(Y,E,F):=\bigl(Y, \tsum_{i\in I} E_i, \tsum_{j\in J} a_j F_j\bigr)
\eqno{(\ref{indestup}.1)}
$$
where $Y$ is a variety  over $\c$, $E_i, F_j$ are
codimension 1 subvarieties and $a_j\in \n$.
(The construction (\ref{summary.k-exmp}) produces a triple
$\bigl(Y, E:=Z, F:=\emptyset\bigr) $. The role of the $F_j$ is to keep track
of the exceptional divisors as we resolve the singularities of $Y$.)

We assume that   $E$ is a simple normal crossing variety and
for every point $p\in E$
there is a  (Euclidean) open neighborhood $p\in Y_p\subset Y$,
 an embedding $\sigma_p : Y_p\into \c^{\dim Y+1}$
whose image can be described as follows.

There are
subsets $I_p\subset I$ and $J_p\subset J$, a natural number $m_p\in \n$
 and   coordinates in $\c^{\dim Y+1}$ called 
$$
\{x_i: i\in I_p\},\ \{y_{rs}:1\leq r, s\leq m_p\},\  \{z_j: j\in J_p\}
\qtq{and} t
$$
(plus possibly other unnamed coordinates)
such that  
 $\sigma_p(Y_p)\subset \c^{\dim Y+1}$ is an open subset of the hypersurface
$$
 \tprod_{i\in I_p}x_i=
t\cdot \det\bigl(y_{rs}:1\leq r, s\leq m_p\bigr)\cdot \tprod_{j\in J_p}z_j^{a_j}.
\eqno{(\ref{indestup}.2)}
$$
Furthermore,   
$$
\begin{array}{l}
\sigma_p(E_i)=(t=x_i=0)\cap \sigma_p(Y_p) \qtq{for} i\in I_p
 \qtq{and}\\
\sigma_p(F_j)=(z_j=0)\cap \sigma_p(Y_p)\qtq{for}j\in J_p.
\end{array}
$$
We do not impose any compatibility condition between the
local equations on overlapping charts.

We say that $(Y,E,F)$
is {\it resolved} at $p$ if $Y$ is smooth at $p$. 
\end{say}

The key technical result of this section  is the following.

\begin{prop}\label{main.res.prop} Let 
$(Y,E,F)$
be a triple as above.
Then there is a resolution of singularities
$\pi: \bigl(Y', E', F'\bigr)\to 
\bigl(Y, E, F\bigr)$
such that
\begin{enumerate}
\item $Y'$ is smooth and $E'$ is a 
simple normal crossing divisor,
\item $ E'=\pi^{-1}(E)$,
\item every stratum of $E'$ is mapped
birationally to a stratum of $E$ and
\item $\pi$ induces an identification
${\mathcal D}(E')={\mathcal D}(E)$.
\end{enumerate}
\end{prop}

Proof.  The resolution will be a composite of explicit blow-ups
of smooth subvarieties (except at the last step).
We use the local equations to describe the blow-up centers locally.
Thus we need to know which 
locally defined subvarieties make sense globally. 
For example,  choosing a
divisor $F_{j_1}$ specifies the local divisor $(z_{j_1}=0)$ at every point
$p\in F_{j_1}$.
Similarly, choosing two divisors $E_{i_1}, E_{i_2}$ gives
the local subvarieties   $(t=x_{i_1}=x_{i_2}=0)$  at every point
$p\in E_{i_1}\cap E_{i_2}$. (Here it is quite important that
the divisors $E_i$ are themselves smooth. The algorithm does not
seem to work if the $E_i$ have self-intersections.) Note that by contrast
$(x_{i_1}=x_{i_2}=0)\subset Y$ defines a local divisor 
which has no global meaning.
Similarly, the vanishing of any of the 
coordinate functions  $y_{rs}$ has no global meaning.

To a point $p\in \sing E$ we associate the local invariant
$$
\mdeg(p):=\bigl(\deg_x(p), \deg_y(p), \deg_z(p)\bigr)=
\bigl(|I_p|, m_p, \tsum_{j\in J_p} a_j\bigr).
$$
It is clear that $\deg_x(p)$ and $ \deg_z(p) $
do not depend on the local coordinates chosen. 
We see in (\ref{detres.say}) that $\deg_y(p)$ is also well defined
if $p\in \sing E$.
The degrees $\deg_x(p), \deg_y(p), \deg_z(p) $
are constructible and upper semi continuous functions on $\sing E$.

Note that $Y$ is smooth at $p$ iff either
$\mdeg (p)=(1,*,*)$ or $\mdeg (p)=(*,0,0)$. 
If $\deg_x(p)=1$ then we can rewrite the equation (\ref{indestup}.2)
as
$$
x'= t\cdot \tprod_{j}z_j^{a_j}
\qtq{where} 
x':=x_1+t\cdot \bigl(1-\det(y_{rs})\bigr)\cdot \tprod_{j}z_j^{a_j},
$$
so if $Y$is smooth then $\bigl(Y, E+ F\bigr)$
has only simple normal crossings along $E$. 
Thus the resolution constructed in Theorem \ref{main.thm.3}
is a log resolution.

The usual method of Hironaka would start by blowing up the {\em highest}
multiplicity points. This introduces new and rather complicated
exceptional divisors and I have not been able to understand
 how the  dual complex changes.

In our case, it turns out to be   better to
look at a  locus where $\deg_y(p)$ is maximal
but instead of maximizing $\deg_x(p)$ or $ \deg_z(p) $
we maximize the dimension. Thus we blow up 
subvarieties along which $Y$ is not equimultiple.
Usually this leads to a morass, but our equations separate
the variables into distinct groups which makes these blow-ups
easy to compute.

One can think of this as  mixing  the main step of the   Hironaka method
with the   order reduction for monomial ideals
(see, for instance, \cite[Step 3 of 3.111]{k-res}).

After some preliminary remarks about
blow-ups of simple normal crossing varieties 
 the proof of (\ref{main.res.prop})
is carried out in a series of steps (\ref{detres.say}--\ref{binres.say}). 

We start with the locus where $\deg_y(p)$ is maximal and by a sequence
of blow-ups we eventually achieve that $\deg_y(p)\leq 1$
for every singular point  $p$. This, however, increases $\deg_z$.
Then in 3 similar steps we lower the maximum of $\deg_z$
until we achieve that $\deg_z(p)\leq 1$ for every singular point   $p$. 
Finally we take care of the  singular points where 
$\deg_y(p)+\deg_z(p)\geq 1$. \qed

\begin{say}[Blowing up simple normal crossing varieties]
\label{snc.blow-up.say}
Let $Z$ be a simple normal crossing variety and
$W\subset Z$ a  subvariety. 
We say that $W$ has {\it simple normal crossing} with $Z$ if
for each point $p\in Z$ there is an open neighborhood $Z_p$,
an embedding $Z_p\into \c^{n+1}$ and subsets
 $I_p, J_p\subset \{0, \dots, n\}$
such that
$$
Z_p=\bigl(\tprod_{i\in I_p} x_i=0\bigr)\qtq{and} 
W\cap Z_p=\bigl(x_j=0: j\in J_p\bigr).
$$
This implies that  
for every stratum $Z_{J}\subset Z$ the intersection
$W\cap Z_J$ is smooth (even scheme theoretically).  

If $W$ has  simple normal crossing with $Z$
then the blow-up $B_WZ$ is again a simple normal crossing variety.
If $W$ is one of the strata of $Z$, then 
${\mathcal D}(B_WZ)$ is obtained from ${\mathcal D}(Z)$
by removing the cell corresponding to $W$ and every other cell whose
closure contains it. Otherwise
${\mathcal D}(B_WZ)={\mathcal D}(Z)$.
(In the terminology of \cite[Sec.2.4]{kk-singbook},
  $B_WZ\to Z$ is a thrifty modification.)

As an example, let $Z=(x_1x_2x_3=0)\subset \c^3$. There are 7 strata
and  ${\mathcal D}(Z)$ is the 2-simplex whose vertices correspond to
the planes $(x_i=0)$.

Let us blow up a point $W=\{p\}\subset Z$ to
 get $B_pZ\subset B_p\c^3$. Note that the exceptional divisor 
$E\subset B_p\c^3$ is {\em not} a part of $B_pZ$
and $B_pZ$ still has 3 irreducible components.

If $p$ is the origin, then the triple intersection is removed
and  ${\mathcal D}(B_pZ)$ is the boundary of the  2-simplex.

If $p$ is not the origin, then $B_pZ$ still has 7 strata
naturally corresponding to the strata of $Z$ and
  ${\mathcal D}(B_pZ)$ is  the  2-simplex.

We will be interested in  situations where $Y$ is a hypersurface in $\c^{n+2}$
and  $Z\subset Y$ is a Cartier divisor that is a 
simple normal crossing variety.
Let $W\subset Y$ be a smooth, irreducible
 subvariety, not contained in $Z$  such that
\begin{enumerate}
\item the scheme theoretic intersection
$W\cap Z$ has simple normal crossing with $Z$
\item $\mult_{Z\cap W}Z=\mult_WY$. (Note that
this holds if $W\subset \sing Y$ and
$\mult_{Z\cap W}Z=2$.)
\end{enumerate}
Choose local coordinates  $(x_0,\dots, x_n,t)$ such that
$W=(x_0=\cdots x_i=0)$ and $Z=(t=0)\subset Y$.
Let $f(x_0,\dots, x_n,t)=0$ be the local equation of $Y$.

Blow up $W$ to get $\pi:B_WY\to Y$. 
Up to permuting the indices $0,\dots, i$, the blow-up $B_WY$ is covered
by coordinate charts described  by
the coordinate change
$$
\bigl(x_0, x_1, \dots, x_i, x_{i+1}, \dots, x_n,t\bigr)=
\big(x'_0, x'_1x'_0, \dots, x'_ix'_0, x_{i+1}, \dots, x_n,t\bigr).
$$
If $\mult_WY=d$ then the local equation of $B_WY$ in the above chart becomes
$$
(x'_0)^{-d}f\big(x'_0, x'_1x'_0, \dots, x'_ix'_0, x_{i+1}, \dots, x_n,t\bigr)=0.
$$
By assumption (2), 
$(x'_0)^{d} $ is also the largest power that divides
$$
f\big(x'_0, x'_1x'_0, \dots, x'_ix'_0, x_{i+1}, \dots, x_n,0\bigr),
$$
hence
$\pi^{-1}(Z)=B_{W\cap Z}Z$.

Observe finally that the conditions (1--2) can not be fulfilled in any
interesting way if $Y$ is smooth. Since we want $Z\cap W$ to be scheme
theoretically smooth,  if $Y$ is smooth then condition (1)
implies that $Z\cap W$ is disjoint from $\sing Z$.

(As an example, let $Y=\c^3$ and $Z=(xyz=0)$. Take $W:=(x=y=z)$. 
Note that $W$ is transversal to every irreducible component of $Z$
but $W\cap Z$ is a non-reduced point. The preimage of
$Z$ in $B_WY$ does not have  simple normal crossings.)

There are, however, plenty of  examples where $Y$ is singular
along $Z\cap W$ and these are exactly the singular points that
we want to resolve.

\end{say}

\begin{say}[Resolving the determinantal part]\label{detres.say}
Let $m$ be the largest size of a determinant occurring at a non-resolved point.
Assume that $m\geq 2$ and let
 $p\in Y$ be a non-resolved point with $m_p=m$. 

Away from $E \cup F $ 
the local equation of $Y$ is 
$$
 \tprod_{i\in I_p}x_i= \det\bigl(y_{rs}:1\leq r, s\leq m\bigr).
$$
Thus,
the singular set of 
$ Y_p\setminus (E \cup F)$ is
$$
\textstyle{\bigcup}_{(i, i')}
\bigl(\rank (y_{rs})\leq m-2\bigr)\cap 
 \bigl( x_{i}=x_{i'}=0\bigr)
$$
where the union runs through all 2-element subsets
$\{i,i'\}\subset I_p$. 
Thus the irreducible components of
$\sing Y\setminus (E \cup F)$
are in natural one-to-one correspondence with the
irreducible components of $\sing E$ and
the value of $m=\deg_y(p)$
is determined by the multiplicity of any of these
irreducible components at $p$.

Pick $i_1, i_2\in I$ and
we work locally with a subvariety 
$$
W'_p(i_1, i_2)
:=\bigl(\rank (y_{rs})\leq m-2\bigr)\cap\bigl(x_{i_1}=x_{i_2}=0\bigr).
$$
Note that $W'_p(i_1, i_2)$ is singular if $m>2$ and
the subset of its highest multiplicity points  is given by $\rank (y_{rs}) =0$.
Therefore the locally defined subvarieties
$$
W_p(i_1, i_2):=\bigl(y_{rs}=0 : 1\leq r,s\leq m\bigr)
\cap\bigl(x_{i_1}=x_{i_2}=0\bigr).
$$
glue together to  a well defined
global smooth subvariety $W:=W(i_1, i_2) $.

$E$ is defined by $(t=0)$ thus $E\cap W $ has the same
local equations as $W_p(i_1, i_2) $. In particular,
$E\cap W $ has simple normal crossings with $E$ and
$E\cap W $ is not a stratum of $E$;
its codimension in the stratum $(x_{i_1}=x_{i_2}=0) $ is $m^2$.

Furthermore, $E$ has multiplicity 2 along $E\cap W $,
hence (\ref{snc.blow-up.say}.2) also holds and so
$$
{\mathcal D}\bigl(B_{E\cap W}\bigr)= {\mathcal D}(E).
$$

We blow up $W \subset Y$. 
We will check that the new triple is again of the form
(\ref{indestup}).
The local degree $\mdeg(p)$
is unchanged over $Y\setminus  W $. The key assertion is that,
over $W $, the maximum value of $\mdeg(p)$
(with respect to the lexicographic ordering) decreases.
By repeating this procedure 
for every irreducible components of $\sing E$,
we  decrease the maximum value of
$\mdeg(p)$.
We can repeat this until we reach $\deg_y(p)\leq 1$ for every 
non-resolved point $p\in Y$.

(Note that this procedure requires an actual ordering of the
irreducible components of $\sing E$, which is a non-canonical choice.
If a finite groups acts on $Y$, 
our resolution  usually can not be chosen  equivariant.)

Now to the local computation of the blow-up.
Fix a  point $p\in W $ and
set $I^*_p:=I_p\setminus\{i_1, i_2\}$.  
We write the local equation of $Y$  as
$$
x_{i_1}x_{i_2}\cdot L=t\cdot \det(y_{rs})\cdot R\qtq{where}
L:=\tprod_{i\in I^*_p}x_i\qtq{and}
R:= \tprod_{j\in J_p}z_j^{a_j}.
$$
Since $W=\bigl(x_{i_1}=x_{i_2}=y_{rs}=0: 1\leq r,s\leq m\bigr)$
there are two types of local charts on  the blow-up.

\begin{enumerate}
\item  There are two  charts of the first type. Up to interchanging
the subscripts $1,2$, these are given by the coordinate change
$$
(x_{i_1},x_{i_2},y_{rs} : 1\leq r,s\leq m)=
(x'_{i_1}, x'_{i_2}x'_{i_1}, y'_{rs}x'_{i_1} : 1\leq r,s\leq m).
$$ 
After setting $z_{w}:=x'_{i_1}$ the new local equation is
$$
 x'_{i_2}\cdot L=t\cdot  \det(y'_{rs})\cdot \bigl(z_{w}^{m^2-2} \cdot R\bigr).
$$  
The exceptional divisor is added to the $F$-divisors with coefficient
$m^2-2$
and the new degree is
$\bigl(\deg_x(p)-1, \deg_y(p), \deg_z(p)+m^2-2\bigr)$.

\item  There are $m^2$  charts of the second type. Up to 
re-indexing the  $m^2$ pairs $(r,s)$  these are given by the coordinate change
$$
(x_{i_1},x_{i_2},y_{rs} : 1\leq r,s\leq m)=
(x'_{i_1}y''_{mm},x'_{i_2}y''_{mm},y'_{rs}y''_{mm} : 1\leq r,s\leq m)
$$
except when $r=s=m$ where we set  $y_{mm}=y''_{mm} $.
It is convenient to set $y'_{mm} =1$ and
$z_w:=y''_{mm} $. Then the new local equation is
$$
x'_{i_1}x'_{i_2}\cdot L=t\cdot  \det\bigl(y'_{rs} : 1\leq r,s\leq m\bigr)
\cdot  \bigl(z_{w}^{m^2-2} \cdot R\bigr).
$$  
Note that  the $(m, m)$ entry of $(y'_{rs})$ is 1. 
By row and column operations we see that
$$
\det \bigl(y'_{rs}: 1\leq r,s,\leq m\bigr)=
\det \bigl(y'_{rs}-y'_{rm}y'_{ms}: 1\leq r,s,\leq m-1\bigr).
$$
By setting $y''_{rs}:=y'_{rs}-y'_{rm}y'_{ms}$ we have new local
equations
$$
x'_{i_1}x'_{i_2}L=t\cdot  
\det \bigl(y''_{rs}: 1\leq r,s,\leq m-1\bigr)
\cdot \bigl(z_{w}^{m^2-2}  \cdot R\bigr)
$$  
 and the new degree is
$\bigl(\deg_x(p), \deg_y(p)-1, \deg_z(p)+m^2-2\bigr)$.

\end{enumerate}
\medskip

{\it Outcome.} After these blow ups we have   a triple
$(Y,E,F)$
such that at non-resolved points the local equations are
$$
 \tprod_{i\in I_p}x_i=
t\cdot y\cdot \tprod_{j\in J_p}z_j^{a_j}
\qtq{or} 
 \tprod_{i\in I_p}x_i=
t\cdot \tprod_{j\in J_p}z_j^{a_j}.
\eqno{(\ref{detres.say}.3)}
$$
(Note that we can not just declare that $y$ is also a $z$-variable.
The $z_j$ are local equations of the divisors $F_j$ while
$(y=0)$ has no global meaning.)
\end{say}

\begin{say}[Resolving the monomial part]\label{monres.say}
Following  (\ref{detres.say}.3), the local equations are
$$
 \tprod_{i\in I_p}x_i=
t\cdot y^c\cdot \tprod_{j\in J_p}z_j^{a_j}\qtq{where $c\in \{0,1\}$.  }
$$
 We lower the degree 
of the $z$-monomial in 3 steps. 

{\it Step} 1. Assume that there is a non-resolved point with
$a_{j_1}\geq 2$. 

The singular set of $F_{j_1}$ is then
$$
\textstyle{\bigcup}_{(i, i')} \bigl(z_{j_1}=x_{i}=x_{i'}=0\bigr)
$$
where the union runs through all 2-element subsets
$\{i,i'\}\subset I$. 
 Pick an irreducible
component of it, call it 
$W(i_1, i_2, j_1):=\bigl(z_{j_1}=x_{i_1}=x_{i_2}=0\bigr) $.

Set $I^*_p:=I_p\setminus\{i_1, i_2\}$, $J^*_p:=J_p\setminus\{j_1\}$ and
 write the local equations as
$$
x_{i_1}x_{i_2}\cdot L=t z_j^{a_j}\cdot R\qtq{where}
L:=\tprod_{i\in I^*_p}x_i\qtq{and}
R:=y^c\cdot \tprod_{j\in J^*_p}z_j^{a_j}.
$$
There are 3 local charts on the blow-up:
\begin{enumerate}
\item $(x_{i_1},x_{i_2},z_j)=(x'_{i_1}, x'_{i_2}x'_{i_1}, z'_jx'_{i_1})$ and,
after setting $z_w:=x'_{i_1}$ 
the new local equation is
$$
 x'_{i_2}\cdot L=t\cdot  z_{w}^{a_j-2}{z'_j}^{a_j}\cdot R.
$$
The new degree is
$\bigl(\deg_x(p)-1, \deg_y(p), \deg_z(p)+a_j-2\bigr)$.
\item Same as above with the subscripts $1,2$ interchanged.
\item $(x_{i_1},x_{i_2},z_j)=(x'_{i_1}z'_j, x'_{i_2}z'_j, z'_j)$
with new local equation
$$
x'_{i_1}x'_{i_2}\cdot L=t \cdot {z'_j}^{a_j-2}\cdot R.
$$
The new degree is
$\bigl(\deg_x(p), \deg_y(p), \deg_z(p)-2\bigr)$.
\end{enumerate}

{\it Step} 2. Assume that there is a non-resolved point with
$a_{j_1}=a_{j_2}=1$. 

The singular set of $F_{j_1}\cap F_{j_2}$ is then
$$
\textstyle{\bigcup}_{(i, i')} \bigl(z_{j_1}=z_{j_2}=x_{i}=x_{i'}=0\bigr).
$$
where the union runs through all 2-element subsets
$\{i,i'\}\subset I$. 
 Pick an irreducible
component of it, call it 
$W(i_1, i_2, j_1, j_2):=\bigl(z_{j_1}=z_{j_2}=x_{i_1}=x_{i_2}=0\bigr) $.

Set $I^*_p:=I_p\setminus\{i_1, i_2\}$, $J^*_p:=J_p\setminus\{j_1, j_2\}$ and
we write the local equations as
$$
x_{i_1}x_{i_2}\cdot L=t z_{j_1}z_{j_2}\cdot R\qtq{where}
L:=\tprod_{i\in I^*_p}x_i\qtq{and}
R:=y^c\cdot \tprod_{j\in J^*_p}z_j^{a_j}.
$$
There are two types of local charts on the blow-up.

\begin{enumerate}
\item In the chart $(x_{i_1},x_{i_2},z_{j_1},z_{j_2})=
(x'_{i_1}, x'_{i_2}x'_{i_1}, z'_{j_1}x'_{i_1},z'_{j_2}x'_{i_1})$
 the new local equation is 
$$
 x'_{i_2}\cdot L=t\cdot  z'_{j_1}z'_{j_2}\cdot R.
$$ 
and the new degree is
$\bigl(\deg_x(p)-1, \deg_y(p), \deg_z(p)\bigr)$.
A similar chart is obtained by interchanging
the subscripts $i_1,i_2$.

\item In the chart $(x_{i_1},x_{i_2},z_{j_1},z_{j_2})=
(x'_{i_1}z'_{j_1}, x'_{i_2}z'_{j_1},z'_{j_1} ,z'_{j_2}z'_{j_1})$.
 the new local equation is 
$$
x'_{i_1}x'_{i_2}\cdot L=t \cdot z'_{j_2} \cdot R.
$$
The new degree is
$\bigl(\deg_x(p), \deg_y(p), \deg_z(p)-1\bigr)$.

A similar chart is obtained by interchanging
the subscripts $j_1,j_2$.
\end{enumerate}

By repeated application of these two steps we are reduced to the case
where $\deg_z(p)\leq 1$  at all  non-resolved points.

{\it Step} 3. Assume that there is a non-resolved point with
$\deg_y(p)=\deg_z(p)=1$. 

The singular set of $Y$ is 
$$
\textstyle{\bigcup}_{(i, i')} \bigl(y=z=x_{i}=x_{i'}=0\bigr).
$$
 Pick an irreducible
component of it, call it 
$W(i_1, i_2):=\bigl(y=z=x_{i_1}=x_{i_2}=0\bigr) $.
The blow up computation is the same as in Step 2.
\medskip

As before we see that at each step the conditions
(\ref{snc.blow-up.say}.1--2) hold, hence
${\mathcal D}(E)$ is unchanged. 
\medskip

{\it Outcome.} After these blow-ups we have   a triple
$(Y,E,F)$
such that at non-resolved points the local equations are
$$
 \tprod_{i\in I_p}x_i= t\cdot y, \quad 
 \tprod_{i\in I_p}x_i= t\cdot z_1 \qtq{or} 
 \tprod_{i\in I_p}x_i= t.
\eqno{(\ref{monres.say}.4)}
$$
As before, the $y$ and $z$ variables have different meaning,
but we can rename $z_1$ as $y$. Thus we have only one 
non-resolved local form left:
$\tprod x_i= t y $.
\end{say}

\begin{say}[Resolving the multiplicity 2 part]\label{binres.say}
 
Here we have a local equation   $x_{i_1}\cdots x_{i_d}=ty$
where $d\geq 2$. We would like to blow up
$(x_{i_1}=y=0)$, but, as we noted, 
this subvariety is not globally defined.
However,  a rare occurrence helps us out.
Usually the blow-up of a smooth subvariety
 determines its center uniquely. However, this is not the case
for codimension 1 centers. Thus we could get a globally well defined
blow-up even from centers that are not globally well defined.

Note that the inverse of $(x_{i_1}=y=0)$ in the local Picard 
group of $Y$ is $E_{i_1}=(x_{i_1}=t=0)$, which is globally defined.
Thus
$$
\proj_Y \tsum_{m\geq 0} \o_Y(mE_{i_1})
$$
is well defined, and locally it is isomorphic to the blow-up
 $B_{(x_{i_1}=y=0)}Y$. (A priori, we would need to take the
normalization of $B_{(x_{i_1}=y=0)}Y$, but it is actually normal.)
Thus we have 2 local charts.

\begin{enumerate}
\item $(x_{i_1},y)=(x'_{i_1}, y'x'_{i_1})$ and
the new local equation is
$\bigl(x_{i_2}\cdots x_{i_d}=ty'\bigr)$.
The new local degree is
$(d-1, 1,0)$.
\item $(x_{i_1},y)=(x'_{i_1}y', y')$ and
the new local equation is
$\bigl(x'_{i_1}\cdot x_{i_2}\cdots x_{i_d}=t\bigr)$.
The new local degree is
$(d, 0,0)$.
\end{enumerate}
\medskip

{\it Outcome.} After all these blow-ups we have   a triple
$\bigl(Y, \tsum_{i\in I} E_i, \tsum_{j\in J} a_j F_j\bigr)$
where  $\tsum_{i\in I} E_i $ is a simple normal crossing divisor and
$Y$ is smooth along   $\tsum_{i\in I} E_i $.
\medskip

This completes the proof of Proposition \ref{main.res.prop}. \qed
\end{say}

\begin{say}[Proof of Theorem \ref{main.thm.2}]
Assume that  $T$ is $\q$-acyclic.
Then, by (\ref{k-fg.thm}) there is a simple normal crossing variety
$Z_T$ such that $H^i\bigl(Z_T, \o_{Z_T}\bigr)=0$ for $i>0$.
Then \cite[Prop.9]{k-t-lc-exmp} shows that,
for $L$ sufficiently ample, 
 the 
singularity $(0\in X_T)$ constructed in (\ref{summary.k-exmp}) 
and (\ref{k-esmp.prop})
is rational. By (\ref{main.res.prop})
we conclude that $\dres(0\in X_T)\cong {\mathcal D}(Z_T)$
is homotopy equivalent to $T$.
\end{say}

\section{Cohen--Macaulay singularities}\label{CM.sec}

\begin{defn}\label{CM.defn}
 {\it Cohen--Macaulay singularities} form the largest
class where Serre duality holds. That is, 
if $X$ is a projective variety of pure dimension $n$ 
then $X$ has  Cohen--Macaulay singularities
iff  $H^i(X, L)$ is dual to $H^{n-i}(X, \omega_X\otimes L^{-1})$
for every line bundle $L$. 
A pleasant property is that if $D\subset X$ is  Cartier divisor
in a scheme then $D$ is Cohen--Macaulay iff $X$ is Cohen--Macaulay
in a neighborhood of $D$.
See 
\cite[pp.184--186]{hartsh} or \cite[Sec.5.5]{km-book} for details.

For local questions it is more convenient to use a characterization
using local cohomology due to
\cite[Sec.3.3]{gro-loc-coh-MR0224620}: $X$ is Cohen--Macaulay  iff
$H^i_x(X, \o_X)=0$ for every $x\in X$ and $i<\dim X$.

Every normal surface is Cohen--Macaulay, so the
topology of the links of Cohen--Macaulay singularities starts to become
 interesting  when $\dim X\geq 3$.
\end{defn}

\begin{defn}\label{perfect.1}
Recall that a group $G$ is called {\it perfect} if 
it has no nontrivial abelian quotients. Equivalently, if
$G=[G,G]$ or if $H_1(G,\z)=0$.

We say that  $G$ is  {\it $\q$-perfect} if every
abelian quotient is torsion. Equivalently, if $H_1(G,\q)=0$.
\end{defn}

The following theorem describes  the  fundamental group of the link of 
Cohen--Macaulay singularities.
Note, however, that the most natural part 
is the equivalence (\ref{cm.main.thm}.1)  $\Leftrightarrow$  
(\ref{cm.main.thm}.5), relating 
the fundamental group of the link 
to the vanishing of  $R^1f_*\o_Y$ for a resolution
$f:Y\to X$. 

\begin{thm} \label{cm.main.thm}
 For a finitely presented group $G$ the following are equivalent.
\begin{enumerate}
\item $G$ is $\q$-perfect (\ref{perfect.1}).
\item $G$ is the fundamental group of the link of an isolated 
Cohen--Macaulay singularity of dimension $= 3$.
\item $G$ is the fundamental group of the link of an isolated 
Cohen--Macaulay singularity of dimension $\geq 3$.
\item $G$ is the fundamental group of the link of a
Cohen--Macaulay singularity whose singular set has codimension $\geq 3$.
\item $G$ is the fundamental group of the link of a 1-rational
singularity (\ref{rtl.defn}).
\end{enumerate}
\end{thm}

Proof. It is clear that (2) $\Rightarrow$ (3)  $\Rightarrow$ (4)
and (\ref{CM=>.1-rat.lem}) shows that (4)  $\Rightarrow$ (5).

The implication 
(5)  $\Rightarrow$  (1) is proved in (\ref{1-rat.link.h1=0.lem}).

Let us prove (1)  $\Rightarrow$  (2).
By (\ref{k-fg.thm}) there is a simple normal crossing variety $Z$ such that
$\pi_1(Z)\cong G$. By a singular version of the Lefschetz
hyperplane theorem (see, for instance,  \cite[Sec.II.1.2]{gm-book}),
by taking  general hyperplane sections we obtain a 
 simple normal crossing surface $S$ such that
$\pi_1(S)\cong G$. Thus  $H^1(S,\q)=0$ and by Hodge theory this implies that 
 $H^1(S,\o_{S})=0$.

By (\ref{k-esmp.prop}) there is a 3--dimensional isolated singularity
$(x\in X)$ with a partial resolution $f:Y\to X$ whose exceptional
divisor is $E\cong S$ and  $R^1f_*\o_Y\cong H^1(E, \o_E)=0$.
In this case the singularities of $Y$ are the simplest possible:
we have only ordinary nodes with equation $(x_1x_2=ty_{11})$.
These are resolved in 1 step by blowing up
$(x_1=t=0)$ and they have no effect on our computations.

Thus $X$ is Cohen--Macaulay by (\ref{CM=i-rat.lem}). \qed

\begin{lem}\label{CM=>.1-rat.lem}
 Let $X$ be  a  normal variety with Cohen--Macaulay singularities
($S_3$ would be sufficient)
and $f:Y\to X$ a resolution of singularities. Then
$\supp R^1f_*\o_Y$ has pure codimension 2. Thus if
$\sing X$ has codimension $\geq 3$ then $R^1f_*\o_Y =0$.
\end{lem}

Proof. By localizing at a generic point of $\supp R^1f_*\o_Y$
(or by taking a generic hyperplane section) 
we may assume that $\supp R^1f_*\o_Y=\{x\}$ is a closed point.
Set $E:=f^{-1}(x)$. There is a Leray spectral sequence
$$
H^i_x\bigl(X, R^jf_*\o_X\bigr)\Rightarrow H^{i+j}_E\bigl(Y, \o_Y).
\eqno{(\ref{CM=>.1-rat.lem}.1)}
$$
By a straightforward duality (see, e.g.\ \cite[10.44]{kk-singbook})
$H^{r}_E\bigl(Y, \o_Y)$ is dual to the stalk of
$R^{n-r}f_*\omega_Y$ which is zero for $r<n$ by \cite{Gra-Rie70b}.
Thus (\ref{CM=>.1-rat.lem}.1) gives an exact sequence
$$
H^1_x\bigl(X, \o_X\bigr)\to H^{1}_E\bigl(Y, \o_Y)
\to H^0_x\bigl(X, R^1f_*\o_X\bigr)\to H^2_x\bigl(X, \o_X\bigr).
$$
If $X$ is Cohen--Macaulay and $\dim X\geq 3$ then
$H^1_x\bigl(X, \o_X\bigr)=H^2_x\bigl(X, \o_X\bigr)=0$,
thus 
$$
\bigl(R^1f_*\o_X\bigr)_x\cong H^0_x\bigl(X, R^1f_*\o_X\bigr)\cong
H^{1}_E\bigl(Y, \o_Y)=0.\qed
$$

For isolated singularities, one has the following converse

\begin{lem} \label{CM=i-rat.lem}
Let $(x\in X)$ be a normal, isolated singularity
with a resolution $f:Y\to X$. Then $X$ is Cohen--Macaulay
iff $R^if_*\o_Y=0$ for $0<i<n-1$.
\end{lem}

Proof. The spectral sequence (\ref{CM=>.1-rat.lem}.1) implies that
we have isomorphisms
$$
R^if_*\o_Y\cong H^i_x(X, \o_X) \qtq{for $0<i<n-1$}
$$
and $H^1_x(X, \o_X)=0$ since $X$ is normal. \qed

\begin{lem} \label{1-rat.link.h1=0.lem}
Let $X$ be  a  normal variety with 1-rational singularities (\ref{rtl.defn})
and $x\in X$ a point with link $L:=L(x\in X)$. Then
$H^1(L, \q)=0$.
\end{lem}

Proof. Let $f:Y\to X$ be a resolution such that
$E:=f^{-1}(x)$ is a simple normal crossing divisor.
By \cite[2.14]{steenbrink} the natural maps
$R^if_*\o_Y\to H^i(E, \o_E)$ are surjective, thus
$H^1(E, \o_E)=0$ hence $H^1(E, \q)=0$ by Hodge theory.

Next we prove that $H^1(E, \q)=H^1(L, \q)$. 
Let $x\in N_X \subset X$ be a neighborhood of $x$ such that $\partial N_X=L$
and $N_Y:=f^{-1}(N_X)$ the corresponding neighborhood of $E$
with boundary $\partial N_Y:=L_Y$. Since $L_Y\to L$ has connected fibers,
$H^1(L, \q)\into H^1(L_Y, \q)$ thus it is enough to prove that
$H^1(L_Y, \q)=0$. The exact cohomology  sequence of the pair
$(N_Y, L_Y)$ gives
$$
0=H^1(E, \q)=H^1(N_Y, \q)\to H^1(L_Y, \q)\to 
H^2(N_Y, L_Y, \q)\stackrel{\alpha}{\to} H^2(N_Y, \q)
$$
By  Poincar\'e duality $H^2(N_Y, L_Y, \q)\cong H_{2n-2}(N_Y,\q)$.
Since $N_Y$ retracts to $E$ we see that $H_{2n-2}(N_Y,\q)$
is freely generated by the classes of exceptional divisors
$E=\cup_i E_i$. The map $\alpha$ sends
$\sum m_i[E_i]$ to $c_1\bigl(\o_{N_Y}(\sum m_iE_i)\bigr)$ 
and we need to show that
the latter are nonzero. This follows from the Hodge index theorem.
\qed

\section{Rational singularities}\label{Ratl.sec}

\begin{defn} \label{rtl.defn}
A quasi projective variety $X$ has {\it rational singularities}
if for one (equivalently every) resolution of singularities $p:Y\to X$
and for every algebraic (or holomorphic) vector bundle $F$ on $X$,
the natural maps  $H^i(X, F)\to H^i(Y, p^*F)$ are isomorphisms.
Thus, for purposes of computing cohomology of vector bundles,
$X$ behaves like a smooth variety.
Rational implies  Cohen--Macaulay.
See \cite[Sec.5.1]{km-book} for details.

A more frequently used 
equivalent definition is the following.
$X$ has  rational singularities iff the higher direct images $R^if_*\o_Y$ are
zero for $i>0$ for one (equivalently every) resolution of 
singularities $p:Y\to X$.

We say that $X$ has {\it 1-rational singularities}
if $R^1f_*\o_Y=0$ for one (equivalently every) resolution of 
singularities $p:Y\to X$.
\end{defn}

\begin{say}[Proof of Theorem \ref{main.thm.2}]
Let $p:Y\to X$ be a resolution of singularities such that
$E_x:=p^{-1}(x)$ is a simple normal crossing divisor. 
As we noted in the proof of (\ref{1-rat.link.h1=0.lem}), 
$R^if_*\o_Y\to H^i(E, \o_E)$ is surjective, thus
 $H^i(E, \o_E)=0$ hence  
$H^i\bigl(\dres(x\in X), \q\bigr)=0$ by (\ref{friedman.lem}). 
Thus $\dres(x\in X) $ is $\q$-acyclic.

Conversely, if $T$ is $\q$-acyclic then Theorem \ref{main.thm.1}
constructs a  singularity which is  rational by (\ref{main.thm.1}.3).\qed
\end{say}

Let $L$ be the link of a rational singularity $(x\in X)$.
Since $X$ is Cohen--Macaulay, we know that
$\pi_1(L)$ is  $\q$-perfect (\ref{cm.main.thm}).
It is not known what else can one say about
fundamental groups of  links of rational singularities,
but the fundamental group of the dual complex can be
completely described.

\begin{defn}\label{perfect}
A group $G$ is called  {\it superperfect}  if 
$H_1(G,\z)=H_2(G,\z)=0$; see \cite{Berrick}. 
We say that  $G$ is 
{\it $\q$-superperfect} if 
$H_1(G,\q)=H_2(G,\q)=0$.
Note that every finite group is $\q$-superperfect.
Other examples are the infinite dihedral group
or $\SL(2,\z)$.
\end{defn}

\begin{cor}\label{rtl.suprrperf.thm}\cite[Thm.42]{k-fg}
  Let $(x\in X)$ be a rational singularity. Then
$\pi_1\bigl(\dres(X)\bigr)$ is  $\q$-super\-per\-fect.
Conversely, for every finitely presented, $\q$-superperfect group $G$
there is a 6-dimensional rational singularity $(x\in X)$ 
such that 
$$
\pi_1\bigl(\dres(X)\bigr)=\pi_1\bigl(\res(X)\bigr)=\pi_1\bigl(L(x\in X)\bigr)
\cong G.
$$
\end{cor}

Proof. By a slight variant of the results of
 \cite{Kervaire, Kervaire-Milnor}, for every finitely presented,
$\q$-super\-per\-fect
group $G$ there is a $\q$-acyclic, 5-dimensional manifold (with boundary)
$M$ whose fundamental group is isomorphic to $G$.
Using this $M$ in (\ref{main.thm.2}) 
we get a rational singularity $(x\in X)$  as desired.

Note that just applying the general construction would give 11 dimensional
examples. See \cite[Sec.7]{k-fg} on how to lower 
the dimension to 6.\footnote{A different construction giving 4 and 5 
dimensional examples is in \cite{k-dual2}.} \qed

\section{Questions and problems}\label{sec.open}

\subsection*{Questions about fundamental groups}{\ }

In principle, for any finitely presented
group $G$  one can follow the proof of \cite{k-fg}
and construct links $L$ such that $\pi_1(L)\cong G$. 
However, in almost all cases, the general methods lead to
very complicated examples. It would be useful to start with some
interesting groups and obtain examples that are understandable.
For example, Higman's group  
$$
H=\langle x_i:  x_i [x_i, x_{i+1}], i\in \z/4\z\rangle
$$ 
is perfect, infinite and contains no proper finite index subgroups  
\cite{Higman}.

\begin{prob} \label{higman.constr}
Find an explicit link whose fundamental group is Higman's group.
(It would be especially interesting to find examples that
occur ``naturally'' in algebraic geometry.)
\end{prob}

Note that our results give links with a given fundamental group
but, as far as we can tell, these groups get killed in 
the larger quasi-projective varieties.
(In particular, we do not answer the question
\cite[p.19]{serre-arbres} whether  Higman's group
can be the fundamental group of a smooth variety.)
This leads to the following.

\begin{ques} Let $G$ be a finitely presented group.
Is there a quasi-projective variety $X$ with an isolated singularity
$(x\in X)$ such that $\pi_1\bigl(L(x\in X)\bigr)\cong G$ and
the natural map
$\pi_1\bigl(L(x\in X)\bigr)\to \pi_1\bigl(X\setminus \{x\}\bigr)$
is an injection?
\end{ques}

As Kapovich pointed out, it is not known if every finitely presented
(or finitely generated) group
occurs as a subgroup of the fundamental group of a smooth
projective or quasi-projective variety.

We saw in (\ref{rtl.suprrperf.thm}) that $\q$-superperfect groups are exactly
 those
that occur as  $\pi_1\bigl(\dres(X)\bigr)$ for rational singularities.
Moreover, every  $\q$-superperfect group can be the fundamental group
of a link of a rational singularity. However, there are
rational singularities such that the  fundamental group of their link is
not  $\q$-superperfect. As an example, let $S$ be a fake projective quadric
whose universal cover is  the 2-disc ${\mathbb D}\times {\mathbb D}$
(cf.\ \cite[Ex.X.13.4]{MR97e:14045}).
Let $C(S)$ be a cone over $S$ with link $L(S)$. 
Then 
$$
H^2\bigl(L(S), \q\bigr)\cong H^2\bigl(S, \q\bigr)/\q \cong\q
$$ and
the universal cover of $L$  is an $\r$-bundle over
 ${\mathbb D}\times {\mathbb D}$
hence contractible. Thus 
$$
H^2\bigl(\pi_1(L(S)), \q\bigr)\cong H^2\bigl(L(S), \q\bigr)\cong \q,
$$
so $\pi_1(L(S)) $  is not  $\q$-superperfect. This leads us to the following,
possibly very hard, question.

\begin{prob} Characterize the fundamental groups of links of
rational singularities.
\end{prob}

In this context it is worthwhile to mention the following.

\begin{conj}[Carlson--Toledo]
 The fundamental group of a  smooth projective variety
is not  $\q$-superperfect (unless it is finite).
\end{conj}

More generally, the original conjecture of Carlson and Toledo
asserts that the image  
$$
\im\bigl[H^2\bigl(\pi_1(X), \q\bigr)\to H^2(X,\q)\bigr]
$$ is nonzero and contains a (possibly degenerate)
K\"ahler class, see  \cite[18.16]{shaf-book}.
For a partial solution see \cite{reznikov-kaehler}.
\medskip

Our examples show that 
for every  finitely presented group $G$ 
there is a {\em reducible} simple normal crossing surface $S$  
 such that $\pi_1(S)\cong G$. By \cite{simp}, for every  finitely presented
group $G$  there is a (very singular) {\em irreducible} variety  $Z$  
 such that $\pi_1(Z)\cong G$. 
It is natural to hope to combine these results.
\cite{2012arXiv1201.3129K} proves that 
for every  finitely presented group $G$ 
there is an irreducible surface $S$ with normal crossing and
Whitney umbrella singularities (also called pinch points,
given locally as  $x^2=y^2z$)
 such that $\pi_1(S)\cong G$.

\begin{prob} \cite{simp} 
What can one say about the fundamental groups of
irreducible surfaces with normal crossing singularities?
\end{prob}

Although closely related, the next question should have a
quite different answer.

\begin{prob} What can one say about the fundamental groups of
normal, projective varieties or surfaces?
Are these two classes of groups the same?
\end{prob}

Many of the known restrictions on fundamental groups of
smooth varieties also apply to normal varieties.
For instance, the theory of
Albanese varieties  implies that
the rank of $H_2(X,\q)$ is even for normal, projective varieties $X$. 
Another example is the following. By \cite{MR900825}
any surjection
 $\pi_1(X)\onto \pi_1(C)$ to the  
fundamental group of a curve $C$ of genus $\geq 2$
factors as
$$
\pi_1(X)\stackrel{g_*}{\to} \pi_1(C'){\onto} \pi_1(C)
$$
where $g:X\to C'$ is a  morphism.
(In general there is no morphism $C'\to C$.)

We claim that this also holds for normal varieties $Y$. 
Indeed, let $\pi:Y'\to Y$ be a resolution of singularities.
Any surjection $\pi_1(Y)\onto \pi_1(C)$ induces
$\pi_1(Y')\onto \pi_1(C)$, hence we get a morphism
$g':Y'\to C'$. Let $B\subset Y'$ be an irreducible curve that is
contacted by $\pi$. Then $\pi_1(B)\to \pi_1(Y)$ is trivial
and so is  $\pi_1(B)\to \pi_1(C)$. If $g'|_B:B\to C'$ is not constant
then the induced map $\pi_1(B)\to \pi_1(C')$ has finite index image.
This is impossible since the composite 
$\pi_1(B)\to \pi_1(C')\to \pi_1(C)$ is trivial.
Thus $g'$ descends to $g:Y\to C'$.

For further such results see \cite{MR983460, MR1076513, MR1098610, MR1463177}.
\medskip

Algebraically one can think of the link as the punctured spectrum
of the Henselisation (or completion) of the local ring of $x\in X$.
Although one can not choose a base point, it should be possible
to define an algebraic  fundamental group.
All the examples in  Theorem  \ref{link.thm}  can be realized
on varieties defined over $\q$. Thus they should have an  
algebraic  fundamental group $\pi^{\rm alg}_1\bigl(L (0\in X_{\q})\bigr)$ 
which is an 
extension of the profinite completion  
of $\pi_1\bigl(L (0\in X)\bigr)$
 and of the  Galois group
$\gal\bigl(\bar{\q}/\q\bigr)$. 

\begin{prob}
Define and describe the possible groups
 $\pi^{\rm alg}_1\bigl(L (0\in X_{\q})\bigr)$.
\end{prob}

\subsection*{Questions about the topology of links}{\ }

We saw that the fundamental groups of links can be quite different from
fundamental groups of quasi-projective varieties. However, our results
say very little about the cohomology or other topological properties of links.
It turns out that links have numerous 
restrictive topological properties. I thank J.~Shaneson and L.~Maxim
for bringing many of these to my attention.

\begin{say}[Which manifolds can be links?] \label{links.which}
Let $M$ be a differentiable manifold that is diffeomeorphic to
the link $L$ 
of an isolated complex singularity of dimension $n$. 
Then $M$ satisfies the following.
\medskip

\ref{links.which}.1. $\dim_{\r}M=2n-1$ is odd and $M$ is orientable.
Resolution of singularities shows that $M$ is cobordant to 0.
\medskip

\ref{links.which}.2.  
The  decomposition  $T_X|_L\cong T_L+N_{L,X}$  shows that $T_M$ is 
stably complex.
In particular, its odd integral Stiefel--Whitney classes are zero
\cite{MR0133137}. (More generally, this holds for orientable  real 
hypersurfaces in  complex manifolds.)
\medskip

\ref{links.which}.3. The cohomology groups $H^i(L,\q)$ carry a natural
mixed Hodge structure; see \cite[Sec.6.3]{PetersSteenbrinkBook}
for a detailed treatment and references. Using these,
\cite{MR928298} proves that the cup product
$H^i(L,\q)\times H^i(L,\q)\to H^{i+j}(L,\q)$ is zero if
$i,j<n$ and $i+j\geq n$. 
In particular, the torus ${\mathbb T}^{2n-1}$ can not be a link.
If $X$ is a smooth projective variety then
$X\times \s^1$ can not be a link.
Further results along this direction are in \cite{MR2454366}.
\medskip

\ref{links.which}.4. By \cite[p.548]{MR1102578}, the
components of the Todd--Hirzebruch L-genus of $M$ vanish
above the middle dimension. 
More generally, the purity of the Chern classes and weight considerations
as in (\ref{links.which}.3) show that the $c_i\bigl(T_X|_L\bigr)$
are torsion above the middle dimension. Thus all
Pontryagin classes of $L$ are torsion above the middle dimension.
See also
\cite{MR2376847, MR2376848} for further results on
the topology of singular algebraic varieties which give
restrictions on links as special cases.
\end{say}

There is no reason to believe that this list is complete
and it would be useful to construct many different links
to get some idea of what other restrictions may hold.

Let $(0\in X)\subset (0\in \c^N)$ be an isolated singularity of dimension $n$
and
$L=X\cap \s^{2N-1}(\epsilon)$ its link.
If $X_0:=X$ is smoothable in a family $\{X_t\subset  \c^N\}$ then
$L$ bounds a Stein manifold   $U_t:=X_t\cap {\mathbb B}^{2N}(\epsilon)$
and $U_t$ is homotopic to an $n$-dimensional compact simplicial complex.
This imposes strong restrictions on the topology of smoothable
links; some of these were used in \cite{MR2454366}.
Interestingly, these restrictions use the integral structure
of the cohomology groups. This leads to the following
intriguing possibility.

\begin{ques} Let $L$ be a link of dimension $2n-1$.
Does $L$ bound a $\q$-homology manifold $U$ (of dimension $2n$)
that is $\q$-homotopic to an $n$-dimensional, finite  simplicial complex?
\end{ques}

There is very little evidence to support the above speculation
but it is consistent with known restrictions
on the topology of links and  it would explain many of them.
On the other hand, I was unable to find such $U$
even in some simple cases. For instance, if $(0\in X)$ is a cone over
an Abelian variety (or a product of curves of genus $\geq 2$)
of dimension $\geq 2$ then algebraic deformations of $X$  do not  produce
such a  $U$.

\medskip
Restricting to the cohomology rings,  here are two simple questions.

\begin{ques}[Cohomology of links] 
Is the sequence of Betti numbers of a complex link arbitrary?
Can one describe the possible algebras  $H^*(L,\q)$?
\end{ques}

\begin{ques}[Cohomology of links of  weighted cones] 
We saw in (\ref{seif.Q.homology.say}) that the first Betti number
of the link of a weighted cone (of dimension $>1$) is even. 
One can ask if this is the
only restriction on the Betti numbers
of a complex link of a  weighted cone.
\end{ques}

Philosophically, one of the main results on the topology of
smooth projective varieties, proved in \cite{dgms, MR0646078}, says that
for   simply connected  varieties
the integral cohomology ring 
and the Pontryagin classes
determine the differentiable structure
up to finite ambiguity. 
It is natural to ask what happens for links.

\begin{ques} \label{fomrlity.ques}
To what extent is the diffeomorphism type of a
simply connected link $L$  determined by the cohomology ring
$H^*(L, \z)$ plus some characteristic classes?
\end{ques}

A positive answer to (\ref{fomrlity.ques}) would imply that
general links are indeed very similar to weighted homogeneous links
and to projective varieties.

 \subsection*{Questions about $\dres(0\in X)$}{\ }

The preprint version contained several questions about 
 dual complexes of dlt
pairs; these are corrected and  solved in
 \cite{dkx}.

\subsection*{Embeddings of simple normal crossing varieties}{\ }

In many contexts it has been a difficulty that not every
variety with simple normal crossing singularities can be realized as a
hypersurface in a smooth variety. See for instance
\cite{fujinobook, 2011arXiv1107.5595B, 2011arXiv1109.3205B, kk-singbook}
for such examples and for various partial solutions.

As we discussed in (\ref{disc.main.steps.say.2}),
recent examples of \cite{2012arXiv1206.1994F, 2012arXiv1206.2475F} show that
the answer to the following may be quite complicated.

\begin{ques} Which proper, complex,  simple normal crossing spaces 
can be realized as 
hypersurfaces in a complex manifold?
\end{ques}

\begin{ques} Which projective simple normal crossing varieties can be 
realized as 
hypersurfaces in a smooth projective variety?
\end{ques}

Note that, in principle it  could happen that there is a
projective simple normal crossing variety that can be realized as a
hypersurface in a complex manifold but not
in a smooth projective variety.
\medskip

Let $Y$ be a smooth variety and $D\subset Y$ a compact divisor.
Let $D\subset N\subset Y$ be a  regular neighborhood
with smooth boundary $\partial N$. If $D$ is the exceptional
divisor of a resolution of an isolated  singularity $x\in X$ then
 $\partial N$ is homeomorphic to the link $L(x\in X)$.
It is clear that $D$ and $c_1\bigl(N_{D,X}\bigr)\in H^2(D, \z)$
determine the boundary  $\partial N$, but I found it very hard to
compute concrete examples.

\begin{prob} Find an effective method to compute
the cohomology or the fundamental group of $\partial N$,
at least when $D$ is a simple normal crossing divisor.
\end{prob}


\begin{thebibliography}{DGMS75}

\bibitem[ABC{\etalchar{+}}96]{abckt}
J.~Amor{\'o}s, M.~Burger, K.~Corlette, D.~Kotschick, and D.~Toledo,
  \emph{Fundamental groups of compact {K}\"ahler manifolds}, Mathematical
  Surveys and Monographs, vol.~44, American Mathematical Society, Providence,
  RI, 1996. \MR{1379330 (97d:32037)}

\bibitem[ABW09]{2009arXiv0902.4234A}
D.~{Arapura}, P.~{Bakhtary}, and J.~{W{\l}odarczyk}, \emph{{The combinatorial
  part of the cohomology of a singular variety}}, ArXiv:0902.4234, 2009.

\bibitem[ABW11]{2011arXiv1102.4370A}
\bysame, \emph{{Weights on cohomology, invariants of singularities, and dual
  complexes}}, ArXiv e-prints (2011).

\bibitem[Art70]{artin}
Michael Artin, \emph{Algebraization of formal moduli. {II}. {E}xistence of
  modifications}, Ann. of Math. (2) \textbf{91} (1970), 88--135. \MR{0260747
  (41 \#5370)}

\bibitem[Bar65]{barden}
D.~Barden, \emph{Simply connected five-manifolds}, Ann. of Math. (2)
  \textbf{82} (1965), 365--385. \MR{MR0184241 (32 \#1714)}

\bibitem[Bea96]{MR97e:14045}
Arnaud Beauville, \emph{Complex algebraic surfaces}, second ed., London
  Mathematical Society Student Texts, vol.~34, Cambridge University Press,
  Cambridge, 1996, Translated from the 1978 French original by R. Barlow, with
  assistance from N. I. Shepherd-Barron and M. Reid. \MR{MR97e:14045}

\bibitem[Ber02]{Berrick}
A.~J. Berrick, \emph{A topologist's view of perfect and acyclic groups},
  Invitations to geometry and topology, Oxf. Grad. Texts Math., vol.~7, Oxford
  Univ. Press, Oxford, 2002, pp.~1--28. \MR{1967745 (2004c:20001)}

\bibitem[BG00]{bg00}
Charles~P. Boyer and Krzysztof Galicki, \emph{On {S}asakian-{E}instein
  geometry}, Internat. J. Math. \textbf{11} (2000), no.~7, 873--909.
  \MR{2001k:53081}

\bibitem[BG08]{bg-book}
\bysame, \emph{Sasakian geometry}, Oxford Mathematical Monographs, Oxford
  University Press, Oxford, 2008. \MR{2382957 (2009c:53058)}

\bibitem[BGK05]{bgk}
Charles~P. Boyer, Krzysztof Galicki, and J{\'a}nos Koll{\'a}r, \emph{Einstein
  metrics on spheres}, Ann. of Math. (2) \textbf{162} (2005), no.~1, 557--580.
  \MR{MR2178969 (2006j:53058)}

\bibitem[BGKT05]{bgkt}
Charles~P. Boyer, Krzysztof Galicki, J{\'a}nos Koll{\'a}r, and Evan Thomas,
  \emph{Einstein metrics on exotic spheres in dimensions 7, 11, and 15},
  Experiment. Math. \textbf{14} (2005), no.~1, 59--64. \MR{2146519
  (2006a:53042)}

\bibitem[BM11]{2011arXiv1107.5595B}
Edward Bierstone and Pierre~D. Milman, \emph{{Resolution except for minimal
  singularities I}}, arXiv.org:1107.5595, 2011.

\bibitem[BP11]{2011arXiv1109.3205B}
Edward Bierstone and Franklin~V. Pacheco, \emph{{Resolution of singularities of
  pairs preserving semi-simple normal crossings}}, arXiv.org:1109.3205, 2011.

\bibitem[Bri66]{briesk-exmps}
Egbert Brieskorn, \emph{Beispiele zur {D}ifferentialtopologie von
  {S}ingularit\"aten}, Invent. Math. \textbf{2} (1966), 1--14. \MR{34 \#6788}

\bibitem[Cai61]{Cairns}
Stewart~S. Cairns, \emph{A simple triangulation method for smooth manifolds},
  Bull. Amer. Math. Soc. \textbf{67} (1961), 389--390. \MR{0149491 (26 \#6978)}

\bibitem[Cat91]{MR1098610}
Fabrizio Catanese, \emph{Moduli and classification of irregular {K}aehler
  manifolds (and algebraic varieties) with {A}lbanese general type fibrations},
  Invent. Math. \textbf{104} (1991), no.~2, 263--289. \MR{1098610 (92f:32049)}

\bibitem[Cat96]{MR1463177}
\bysame, \emph{Fundamental groups with few relations}, Higher-dimensional
  complex varieties ({T}rento, 1994), de Gruyter, Berlin, 1996, pp.~163--165.
  \MR{1463177 (98i:32047)}

\bibitem[CMS08a]{MR2376847}
Sylvain~E. Cappell, Laurentiu~G. Maxim, and Julius~L. Shaneson, \emph{Euler
  characteristics of algebraic varieties}, Comm. Pure Appl. Math. \textbf{61}
  (2008), no.~3, 409--421. \MR{2376847 (2009f:14038)}

\bibitem[CMS08b]{MR2376848}
\bysame, \emph{Hodge genera of algebraic varieties. {I}}, Comm. Pure Appl.
  Math. \textbf{61} (2008), no.~3, 422--449. \MR{2376848 (2009f:14039)}

\bibitem[Coh73]{MR0362320}
Marshall~M. Cohen, \emph{A course in simple-homotopy theory}, Springer-Verlag,
  New York, 1973, Graduate Texts in Mathematics, Vol. 10. \MR{0362320 (50
  \#14762)}

\bibitem[Cor92]{Corson}
Jon~Michael Corson, \emph{Complexes of groups}, Proc. London Math. Soc. (3)
  \textbf{65} (1992), no.~1, 199--224. \MR{1162493 (93h:57003)}

\bibitem[CS91]{MR1102578}
Sylvain~E. Cappell and Julius~L. Shaneson, \emph{Stratifiable maps and
  topological invariants}, J. Amer. Math. Soc. \textbf{4} (1991), no.~3,
  521--551. \MR{1102578 (92d:57024)}

\bibitem[CS08]{Corlette-Simpson}
Kevin Corlette and Carlos Simpson, \emph{On the classification of rank-two
  representations of quasiprojective fundamental groups}, Compos. Math.
  \textbf{144} (2008), no.~5, 1271--1331. \MR{2457528 (2010i:14006)}

\bibitem[Dem88]{demazure}
Michel Demazure, \emph{Anneaux gradu\'es normaux}, Introduction \`a la
  th\'eorie des singularit\'es, II, Travaux en Cours, vol.~37, Hermann, Paris,
  1988, pp.~35--68. \MR{91k:14004}

\bibitem[dFKX12]{dkx}
Tommaso de~Fernex, J{\'a}nos Koll{\'a}r, and Chenyang Xu, \emph{{The dual
  complex of singularities}}, ArXiv e-prints (2012).

\bibitem[DGMS75]{dgms}
Pierre Deligne, Phillip Griffiths, John Morgan, and Dennis Sullivan, \emph{Real
  homotopy theory of {K}\"ahler manifolds}, Invent. Math. \textbf{29} (1975),
  no.~3, 245--274. \MR{MR0382702 (52 \#3584)}

\bibitem[DH88]{MR928298}
Alan~H. Durfee and Richard~M. Hain, \emph{Mixed {H}odge structures on the
  homotopy of links}, Math. Ann. \textbf{280} (1988), no.~1, 69--83. \MR{928298
  (89c:14012)}

\bibitem[Dol83]{dol-link}
Igor~V. Dolgachev, \emph{On the link space of a {G}orenstein quasihomogeneous
  surface singularity}, Math. Ann. \textbf{265} (1983), no.~4, 529--540.
  \MR{721886 (85k:32024)}

\bibitem[DPS09]{DPS}
Alexandru Dimca, {S}tefan Papadima, and Alexander~I. Suciu, \emph{Topology and
  geometry of cohomology jump loci}, Duke Math. J. \textbf{148} (2009), no.~3,
  405--457. \MR{2527322 (2011b:14047)}

\bibitem[FM83]{friedman-etal}
Robert Friedman and David~R. Morrison (eds.), \emph{The birational geometry of
  degenerations}, Progr. Math., vol.~29, Birkh\"auser Boston, Mass., 1983.
  \MR{690262 (84g:14032)}

\bibitem[Fuj09]{fujinobook}
Osamu Fujino, \emph{Introduction to the log minimal model program for log
  canonical pairs}, arXiv.org:0907.1506, 2009.

\bibitem[Fuj12a]{2012arXiv1206.1994F}
Kento Fujita, \emph{{Simple normal crossing Fano varieties and log Fano
  manifolds}}, ArXiv e-prints (2012).

\bibitem[Fuj12b]{2012arXiv1206.2475F}
\bysame, \emph{{The Mukai conjecture for log Fano manifolds}}, ArXiv e-prints
  (2012).

\bibitem[FZ03]{fl-za}
Hubert Flenner and Mikhail Zaidenberg, \emph{Normal affine surfaces with {$\Bbb
  C\sp \ast$}-actions}, Osaka J. Math. \textbf{40} (2003), no.~4, 981--1009.
  \MR{2 020 670}

\bibitem[GL91]{MR1076513}
Mark Green and Robert Lazarsfeld, \emph{Higher obstructions to deforming
  cohomology groups of line bundles}, J. Amer. Math. Soc. \textbf{4} (1991),
  no.~1, 87--103. \MR{MR1076513 (92i:32021)}

\bibitem[GM88]{gm-book}
Mark Goresky and Robert MacPherson, \emph{Stratified {M}orse theory},
  Ergebnisse der Mathematik und ihrer Grenzgebiete (3), vol.~14,
  Springer-Verlag, Berlin, 1988. \MR{932724 (90d:57039)}

\bibitem[Gor80]{MR576865}
Gerald~Leonard Gordon, \emph{On a simplicial complex associated to the
  monodromy}, Trans. Amer. Math. Soc. \textbf{261} (1980), no.~1, 93--101.
  \MR{576865 (81j:32017)}

\bibitem[GR70]{Gra-Rie70b}
Hans Grauert and Oswald Riemenschneider, \emph{Verschwindungss\"atze f\"ur
  analytische {K}ohomologiegruppen auf komplexen {R}\"aumen}, Invent. Math.
  \textbf{11} (1970), 263--292. \MR{MR0302938 (46 \#2081)}

\bibitem[Gro67]{gro-loc-coh-MR0224620}
Alexander Grothendieck, \emph{Local cohomology}, Lecture Notes in Mathematics,
  Vol. 41, Springer-Verlag, Berlin, 1967. \MR{0224620 (37 \#219)}

\bibitem[Gro68]{sga2}
\bysame, \emph{Cohomologie locale des faisceaux coh\'erents et th\'eor\`emes de
  {L}efschetz locaux et globaux {$(SGA$} {$2)$}}, North-Holland Publishing Co.,
  Amsterdam, 1968, Augment\'e d'un expos\'e par Mich\`ele Raynaud, S\'eminaire
  de G\'eom\'etrie Alg\'ebrique du Bois-Marie, 1962, Advanced Studies in Pure
  Mathematics, Vol. 2. \MR{0476737 (57 \#16294)}

\bibitem[Gro89]{MR983460}
Michel Gromov, \emph{Sur le groupe fondamental d'une vari\'et\'e
  k\"ahl\'erienne}, C. R. Acad. Sci. Paris S\'er. I Math. \textbf{308} (1989),
  no.~3, 67--70. \MR{983460 (90i:53090)}

\bibitem[GS75]{gri-sch}
Phillip Griffiths and Wilfried Schmid, \emph{Recent developments in {H}odge
  theory: a discussion of techniques and results}, Discrete subgroups of {L}ie
  groups and applicatons to moduli ({I}nternat. {C}olloq., {B}ombay, 1973),
  Oxford Univ. Press, Bombay, 1975, pp.~31--127. \MR{0419850 (54 \#7868)}

\bibitem[Har77]{hartsh}
Robin Hartshorne, \emph{Algebraic geometry}, Springer-Verlag, New York, 1977,
  Graduate Texts in Mathematics, No. 52. \MR{0463157 (57 \#3116)}

\bibitem[Har95]{Harris95}
Joe Harris, \emph{Algebraic geometry}, Graduate Texts in Mathematics, vol. 133,
  Springer-Verlag, New York, 1995, A first course, Corrected reprint of the
  1992 original. \MR{MR1416564 (97e:14001)}

\bibitem[Hat02]{hatcher}
Allen Hatcher, \emph{Algebraic topology}, Cambridge University Press,
  Cambridge, 2002. \MR{1867354 (2002k:55001)}

\bibitem[Hig51]{Higman}
Graham Higman, \emph{A finitely generated infinite simple group}, J. London
  Math. Soc. \textbf{26} (1951), 61--64. \MR{0038348 (12,390c)}

\bibitem[Hir62]{Hirsch}
Morris~W. Hirsch, \emph{Smooth regular neighborhoods}, Ann. of Math. (2)
  \textbf{76} (1962), 524--530. \MR{0149492 (26 \#6979)}

\bibitem[Ish85]{Ishii85}
Shihoko Ishii, \emph{On isolated {G}orenstein singularities}, Math. Ann.
  \textbf{270} (1985), no.~4, 541--554. \MR{MR776171 (86j:32024)}

\bibitem[{Kap}12]{2012arXiv1201.3129K}
M.~{Kapovich}, \emph{{Dirichlet fundamental domains and complex-projective
  varieties}}, ArXiv e-prints (2012).

\bibitem[Ker69]{Kervaire}
Michel~A. Kervaire, \emph{Smooth homology spheres and their fundamental
  groups}, Trans. Amer. Math. Soc. \textbf{144} (1969), 67--72. \MR{0253347 (40
  \#6562)}

\bibitem[KK11]{k-fg}
Michael {Kapovich} and J{\'a}nos Koll{\'a}r, \emph{{Fundamental groups of links
  of isolated singularities}}, Journal AMS (to appear) ArXiv e-prints (2011).

\bibitem[KM63]{Kervaire-Milnor}
Michel~A. Kervaire and John~W. Milnor, \emph{Groups of homotopy spheres. {I}},
  Ann. of Math. (2) \textbf{77} (1963), 504--537. \MR{0148075 (26 \#5584)}

\bibitem[KM98a]{kap-mil}
Michael Kapovich and John~J. Millson, \emph{On representation varieties of
  {A}rtin groups, projective arrangements and the fundamental groups of smooth
  complex algebraic varieties}, Inst. Hautes \'Etudes Sci. Publ. Math. (1998),
  no.~88, 5--95 (1999). \MR{1733326 (2001d:14024)}

\bibitem[KM98b]{km-book}
J{\'a}nos Koll{\'a}r and Shigefumi Mori, \emph{Birational geometry of algebraic
  varieties}, Cambridge Tracts in Mathematics, vol. 134, Cambridge University
  Press, Cambridge, 1998, With the collaboration of C. H. Clemens and A. Corti,
  Translated from the 1998 Japanese original. \MR{1658959 (2000b:14018)}

\bibitem[Kob63]{MR0154235}
Shoshichi Kobayashi, \emph{Topology of positively pinched {K}aehler manifolds},
  T\^ohoku Math. J. (2) \textbf{15} (1963), 121--139. \MR{0154235 (27 \#4185)}

\bibitem[Kol95]{shaf-book}
J{\'a}nos Koll{\'a}r, \emph{Shafarevich maps and automorphic forms}, M. B.
  Porter Lectures, Princeton University Press, Princeton, NJ, 1995. \MR{1341589
  (96i:14016)}

\bibitem[Kol05]{k-em2}
\bysame, \emph{Einstein metrics on five-dimensional {S}eifert bundles}, J.
  Geom. Anal. \textbf{15} (2005), no.~3, 445--476. \MR{MR2190241 (2007c:53056)}

\bibitem[Kol06]{k-circ}
\bysame, \emph{Circle actions on simply connected 5-manifolds}, Topology
  \textbf{45} (2006), no.~3, 643--671. \MR{2218760 (2006m:57044)}

\bibitem[Kol07a]{k-em1}
\bysame, \emph{Einstein metrics on connected sums of {$S\sp 2\times S\sp 3$}},
  J. Differential Geom. \textbf{75} (2007), no.~2, 259--272. \MR{MR2286822
  (2007k:53061)}

\bibitem[Kol07b]{k-res}
\bysame, \emph{Lectures on resolution of singularities}, Annals of Mathematics
  Studies, vol. 166, Princeton University Press, Princeton, NJ, 2007.
  \MR{2289519}

\bibitem[Kol09]{k-em3}
\bysame, \emph{Positive {S}asakian structures on 5-manifolds}, Riemannian
  topology and geometric structures on manifolds, Progr. Math., vol. 271,
  Birkh\"auser Boston, Boston, MA, 2009, pp.~93--117. \MR{2494170
  (2010i:53077)}

\bibitem[Kol11]{k-t-lc-exmp}
\bysame, \emph{New examples of terminal and log canonical singularities},
  arXiv:1107.2864, 2011.

\bibitem[Kol12]{k-q}
\bysame, \emph{Quotients by finite equivalence relations}, Current developments
  in algebraic geometry, Math. Sci. Res. Inst. Publ., vol.~59, Cambridge Univ.
  Press, Cambridge, 2012, With an appendix by Claudiu Raicu, pp.~227--256.
  \MR{2931872}

\bibitem[Kol13a]{k-dual2}
\bysame, \emph{Simple normal crossing varieties with prescribed dual complex},
  ArXiv e-prints (2013).

\bibitem[Kol13b]{kk-singbook}
\bysame, \emph{Singularities of the minimal model program}, Cambridge
  University Press, Cambridge, 2013, With the collaboration of S. Kov\'acs.

\bibitem[Kul77]{MR0506296}
Vik.~S. Kulikov, \emph{Degenerations of {$K3$} surfaces and {E}nriques
  surfaces}, Izv. Akad. Nauk SSSR Ser. Mat. \textbf{41} (1977), no.~5,
  1008--1042, 1199. \MR{0506296 (58 \#22087b)}

\bibitem[Mas61]{MR0133137}
W.~S. Massey, \emph{Obstructions to the existence of almost complex
  structures}, Bull. Amer. Math. Soc. \textbf{67} (1961), 559--564. \MR{0133137
  (24 \#A2971)}

\bibitem[Mor78]{Morgan}
John~W. Morgan, \emph{The algebraic topology of smooth algebraic varieties},
  Inst. Hautes \'Etudes Sci. Publ. Math. (1978), no.~48, 137--204. \MR{516917
  (80e:55020)}

\bibitem[Mum61]{mumf-top}
David Mumford, \emph{The topology of normal singularities of an algebraic
  surface and a criterion for simplicity}, Inst. Hautes \'Etudes Sci. Publ.
  Math. (1961), no.~9, 5--22. \MR{0153682 (27 \#3643)}

\bibitem[Neu81]{neumann-pl}
Walter~D. Neumann, \emph{A calculus for plumbing applied to the topology of
  complex surface singularities and degenerating complex curves}, Trans. Amer.
  Math. Soc. \textbf{268} (1981), no.~2, 299--344. \MR{632532 (84a:32015)}

\bibitem[OW75]{or-wa}
Peter Orlik and Philip Wagreich, \emph{Seifert {$n$}-manifolds}, Invent. Math.
  \textbf{28} (1975), 137--159. \MR{50 \#13596}

\bibitem[Pay09]{payne09}
Sam Payne, \emph{Lecture at {MSRI}},
  http://www.msri.org/web/msri/online-videos/-/video/showVideo/3674, 2009.

\bibitem[{Pay}11]{payne11}
Sam {Payne}, \emph{{Boundary complexes and weight filtrations}}, ArXiv e-prints
  (2011).

\bibitem[Per77]{MR0466149}
Ulf Persson, \emph{On degenerations of algebraic surfaces}, Mem. Amer. Math.
  Soc. \textbf{11} (1977), no.~189, xv+144. \MR{0466149 (57 \#6030)}

\bibitem[Pin77]{pinkham}
H.~Pinkham, \emph{Normal surface singularities with {$C\sp*$} action}, Math.
  Ann. \textbf{227} (1977), no.~2, 183--193. \MR{55 \#5623}

\bibitem[PP08]{MR2454366}
Patrick Popescu-Pampu, \emph{On the cohomology rings of holomorphically
  fillable manifolds}, Singularities {II}, Contemp. Math., vol. 475, Amer.
  Math. Soc., Providence, RI, 2008, pp.~169--188. \MR{2454366 (2010h:32039)}

\bibitem[PS08]{PetersSteenbrinkBook}
Chris A.~M. Peters and Joseph H.~M. Steenbrink, \emph{Mixed {H}odge
  structures}, Ergebnisse der Mathematik und ihrer Grenzgebiete. 3. Folge.,
  vol.~52, Springer-Verlag, Berlin, 2008. \MR{MR2393625}

\bibitem[Rez02]{reznikov-kaehler}
Alexander Reznikov, \emph{The structure of {K}\"ahler groups. {I}. {S}econd
  cohomology}, Motives, polylogarithms and {H}odge theory, {P}art {II}
  ({I}rvine, {CA}, 1998), Int. Press Lect. Ser., vol.~3, Int. Press,
  Somerville, MA, 2002, pp.~717--730. \MR{1978716 (2004c:32042)}

\bibitem[Sco83]{scott}
Peter Scott, \emph{The geometries of {$3$}-manifolds}, Bull. London Math. Soc.
  \textbf{15} (1983), no.~5, 401--487. \MR{84m:57009}

\bibitem[Sei32]{seif}
Herbert Seifert, \emph{Topologie dreidimensionaler gefaserte {R}{\"a}ume}, Acta
  Math. \textbf{60} (1932), 148--238.

\bibitem[Ser77]{serre-arbres}
Jean-Pierre Serre, \emph{Arbres, amalgames, {${\rm SL}_{2}$}}, Soci\'et\'e
  Math\'ematique de France, Paris, 1977, Avec un sommaire anglais,
  R{\'e}dig{\'e} avec la collaboration de Hyman Bass, Ast{\'e}risque, No. 46.
  \MR{0476875 (57 \#16426)}

\bibitem[Sim10]{simp}
Carlos Simpson, \emph{Local systems on proper algebraic {V}-manifolds},
  arXiv1010.3363, 2010.

\bibitem[Siu87]{MR900825}
Yum~Tong Siu, \emph{Strong rigidity for {K}\"ahler manifolds and the
  construction of bounded holomorphic functions}, Discrete groups in geometry
  and analysis ({N}ew {H}aven, {C}onn., 1984), Progr. Math., vol.~67,
  Birkh\"auser Boston, Boston, MA, 1987, pp.~124--151. \MR{900825 (89i:32044)}

\bibitem[Sma62]{smale}
Stephen Smale, \emph{On the structure of {$5$}-manifolds}, Ann. of Math. (2)
  \textbf{75} (1962), 38--46. \MR{25 \#4544}

\bibitem[Ste83]{steenbrink}
J.~H.~M. Steenbrink, \emph{Mixed {H}odge structures associated with isolated
  singularities}, Singularities, {P}art 2 ({A}rcata, {C}alif., 1981), Proc.
  Sympos. Pure Math., vol.~40, Amer. Math. Soc., Providence, RI, 1983,
  pp.~513--536. \MR{713277 (85d:32044)}

\bibitem[Ste08]{MR2399025}
D.~A. Stepanov, \emph{A note on resolution of rational and hypersurface
  singularities}, Proc. Amer. Math. Soc. \textbf{136} (2008), no.~8,
  2647--2654. \MR{2399025 (2009g:32060)}

\bibitem[Sul77]{MR0646078}
Dennis Sullivan, \emph{Infinitesimal computations in topology}, Inst. Hautes
  \'Etudes Sci. Publ. Math. (1977), no.~47, 269--331 (1978). \MR{0646078 (58
  \#31119)}

\bibitem[Thu07]{MR2320738}
Amaury Thuillier, \emph{G\'eom\'etrie toro\"\i dale et g\'eom\'etrie analytique
  non archim\'edienne. {A}pplication au type d'homotopie de certains sch\'emas
  formels}, Manuscripta Math. \textbf{123} (2007), no.~4, 381--451. \MR{2320738
  (2008g:14038)}

\end{thebibliography}

\newcommand{\etalchar}[1]{$^{#1}$}
\def\cprime{$'$} \def\cprime{$'$} \def\cprime{$'$} \def\cprime{$'$}
  \def\cprime{$'$} \def\cprime{$'$} \def\dbar{\leavevmode\hbox to
  0pt{\hskip.2ex \accent"16\hss}d} \def\cprime{$'$} \def\cprime{$'$}
  \def\polhk#1{\setbox0=\hbox{#1}{\ooalign{\hidewidth
  \lower1.5ex\hbox{`}\hidewidth\crcr\unhbox0}}} \def\cprime{$'$}
  \def\cprime{$'$} \def\cprime{$'$} \def\cprime{$'$}
  \def\polhk#1{\setbox0=\hbox{#1}{\ooalign{\hidewidth
  \lower1.5ex\hbox{`}\hidewidth\crcr\unhbox0}}} \def\cdprime{$''$}
  \def\cprime{$'$} \def\cprime{$'$} \def\cprime{$'$} \def\cprime{$'$}
\providecommand{\bysame}{\leavevmode\hbox to3em{\hrulefill}\thinspace}
\providecommand{\MR}{\relax\ifhmode\unskip\space\fi MR }
\providecommand{\MRhref}[2]{%
  \href{http://www.ams.org/mathscinet-getitem?mr=#1}{#2}
}
\providecommand{\href}[2]{#2}

\noindent Princeton University, Princeton NJ 08544-1000

{\begin{verbatim}kollar@math.princeton.edu\end{verbatim}}

\end{document}